\documentclass[10pt, a4paper]{article}
\usepackage[utf8]{inputenc}

\usepackage[T1]{fontenc}
\usepackage{amsmath, amssymb, amsthm}
\newtheorem{pred}{Proposition}[section]
\theoremstyle{remark}
\newtheorem*{remark}{Remark}
\numberwithin{equation}{section}
\usepackage{mathtools}
\usepackage{array}
\usepackage{hyperref}
\usepackage[english]{babel}
\usepackage{booktabs}
\usepackage{latexsym}
\usepackage{bm}
\usepackage[pdftex]{graphicx}
\usepackage{geometry}
\usepackage{breqn}
\usepackage{multicol}
\usepackage[table]{xcolor}

\pdfoutput=1
\usepackage{etoolbox}
\AtBeginEnvironment{tabular}{\linespread{0.8}\selectfont}

\newcommand{\msr}[5][:]{\prescript{#2}{#3}{#1}^{#4}_{#5}{}}
\newcommand{\mstr}[4]{\prescript{#1}{#2}{!}^{#3}_{#4}{}}

\title{The Ultra-Radical: Analytic Continuation, Branching, and Stability of the Principal Branch}
\author{Sergey V. Berezin}
\date{Version 4.0, \today}
\begin{document}

\maketitle

\begin{center}
\textbf{Affiliation}: Municipal Budgetary Institution, Iglino, Republic of Bashkortostan, Russia \\
\textbf{ORCID}: \href{https://orcid.org/0000-0001-8086-8288}{0000-0001-8086-8288} \\
\textbf{Email}: \href{mailto:bkcru@bk.ru}{bkcru@bk.ru} \\
\textbf{License}: Creative Commons Attribution 4.0 International (CC BY 4.0)
\end{center}

\begin{abstract}
We study the ultra-radical $\sqrt[{n;a;b}]{x}$, the multi-valued solution to $y^{a}=1+axy^{b}$. 
Inside the convergence radius $|x|<R$, every branch is given by a Master--J power series; 
for $|x|\ge R$, analytic continuation requires switching to one of two \emph{conjugate series}.

We introduce a \textbf{deterministic geometric criterion} that selects, for each branch index $n$, the correct conjugate series, thereby eliminating heuristic search and guaranteeing branch continuity across $|x|=R$.

\textbf{Key finding:} Only the \emph{principal branch} ($n=0$) remains continuous when the parameters $a$, $b$, and $x$ vary smoothly. This includes the critical limits $a\to0$ (transition to an exponential equation) and $b\to0$ (transition to a binomial root), where the principal branch converges to the corresponding classical solution. In contrast, branches with $n\neq0$ exhibit oscillatory divergence as $a\to0$ and lose their identity in these limits.

This structural continuity singles out the principal branch for applications where parameters may vary with the system's state, such as in nonlinear media with field‑dependent exponents or adaptive dynamical systems.
\end{abstract}

\tableofcontents
\section{Introduction}
Consider the \textbf{master function} $y = M(m; \alpha; \beta; x)$, defined by the series
\begin{equation}
M(m; \alpha; \beta; x) = m + x + \sum_{\ell=2}^{\infty} \frac{x^{\ell}}{\ell!} \prod_{\gamma=1}^{\ell-1} (m - \alpha \gamma + \beta \ell). \label{eq_master_series}
\end{equation}
This is a multi-parameter family of solutions, continuously depending on parameters $\alpha, \beta$ and the independent argument $x$.

When $\alpha = 0$ or $\beta = 0$, the master function is a power series representing solutions of classical elementary or transcendental equations, where $\alpha = a$, $\beta = b$:
\begin{itemize}
    \item For $m = 0$: $y = x$ ($a = b = 0$), $y = x e^{by}$ ($a = 0, b \neq 0$), $y = \frac{\ln(1 + ax)}{a}$ ($a \neq 0, b = 0$).
    \item For $m = 1$: $y = e^{x}$ ($a = b = 0$), $y = e^{xy^{b}}$ ($a = 0, b \neq 0$), $y = (1 + ax)^{1/a}$ ($a \neq 0, b = 0$).
\end{itemize}
For $a \neq 0$ and $b \neq 0$ with $m = 1$, the master function is the power series expansion of the \textbf{ultra-radical}
\begin{equation}
y_n = \sqrt[{n; a; b}]{x} \equiv U_{n}(a;b;x), \qquad n \in \mathbb{Z}, \label{eq_ultra_def_en}
\end{equation}
defined as the multivalued solution of the algebrao-transcendental equation
\begin{equation}
y = (1 + a x y^{b})^{1/a}. \label{eq_ultra_eq_en}
\end{equation}
Its power-series representation possesses a \textbf{finite radius of convergence}
\begin{equation}
R = \left| b^{-\frac{b}{a}} \cdot (b - a)^{\frac{b - a}{a}} \right|, \qquad a \neq 0,\; b \neq 0,\; b \neq a,
\label{eq_radius}
\end{equation}
which for real parameters $a, b \in \mathbb{R}$ simplifies to $R = \frac{|1 - a/b|^{b/a}}{|b - a|}$.

The ultra-radical is always defined via the master series $y = \sqrt[{n; a; b}]{x} = vM(1;\alpha;\beta;X_J V/(\alpha q_J))$,
where $v=e^{f}$, $V=e^{\beta f}$, $f=\frac{\ln|q_J/p_J|+i[\arg(q_J/p_J)+2\pi N]}{\alpha}$. For $|x|<R$, the direct expansion is used,
while for $|x|\ge R$, one of two conjugate representations with different
parameters $\alpha,\beta$ is employed, the choice depending on the branch number $n$:
\[
\begin{array}{c|c|c|c|c|c|c}
\text{Domain} & \text{Param.} & \alpha & \beta & p_J & q_J & X_J \\
\hline
|x|<R\ (\text{inside}) J=1 & N=n\in\mathbb{Z} & a & b & 1 & 1 & ax \\[10pt]
|x|\ge R\ (\text{case h}) J=2 & N=h\in\mathbb{Z} & b-a & -a & -ax & -1 & 1 \\[10pt]
|x|\ge R\ (\text{case k}) J=3 & N=k\in\mathbb{Z} & -b & a-b & -1 & ax & -1
\end{array}
\]
\label{sec_geometric_criterion}
\textbf{Sector centers} and \textbf{candidates} on the imaginary axis:
\[
L(n)=\operatorname{Im}(b f(n)),\quad
L(h)=\operatorname{Im}(b f(h)),\quad  
L(k)=\operatorname{Im}(b f(k)),
\]
For real $a,b$, the expressions simplify to:
\[
L(n)=\frac{2\pi b n}{a},\quad
L(h)=\frac{b[\arg(1/(ax))+2\pi h]}{b-a},\quad
L(k)=\frac{b[\arg(-ax)+2\pi k]}{-b}.
\]

The centers $L(n)$ partition the imaginary axis into sectors $[g_n, G_n]=[L(n)-\pi b/a,\ L(n)+\pi b/a]$.

The \textbf{angular criterion} selects those $h$ or $k$ for which $L(h)$ or $L(k)$ fall into the sector $[g_n, G_n]$ of branch $n$, ensuring continuous analytic continuation.

Note also the connection between the behavior for $|x|\ge R$ and the limit $a\to0$. The two conjugate representations of the principal branch — $h=0$ and $k=0$ — tend, as $a\to0$, to the two distinct real branches of the generalized Lambert function $y=e^{x y^{b}}$ (where they exist).
\textbf{Subsequent developments.} This paper is the first in a series of three works. The algebra of master numbers, the merge operation for series, and Taylor series of ultrafunctions are systematically presented in the second paper~\cite{Berezin2026_arXiv}; a rigorous proof of the sector criterion, its extension to the ultralogarithm, and its generalization to complex exponents are given in the third paper~\cite{Berezin2026_ULog}. For the current statements of the criterion and the boundary convention, the reader should consult these later publications.

\section{Examples of Algorithm Operation}
\noindent
\textbf{Remark on integer real exponents}. Since in the considered examples all exponents $a,b$ are integer real numbers, the corresponding Riemann surface of the ultra-radical has a finite number of sheets, which can be represented as a closed circle. This simplifies visualization and allows clear demonstration of the angular criterion's operation.

\subsection{Example 1: Continuation Arguments Inside Sectors}

Consider the ultra-radical $\sqrt[{n;5;2}]{x}$, solving the equation $y^{5} = 1 + 5xy^{2}$.

\textbf{Problem}: Find asymptotic expansion of branch ($n=2$) for $x=7$.

\textbf{Solution}:
\begin{enumerate}
\item Find sector center: $\arg\left( v(n = 2) \right) = \frac{4\pi}{5} = 144^\circ$
\item Determine sector boundaries: $\left( \frac{4}{5} \pm \frac{2}{5} \right)\pi = [108^\circ, 180^\circ]$
\item Select $y(h)$ or $y(k)$ whose $\arg(v)$ lies within this sector
\end{enumerate}

After calculations:
\begin{align*}
\arg\left( v(h = 0) \right) &= 0^\circ,\\
\arg\left( v(h = 1) \right) &= \frac{2\pi}{- 3} = 300^\circ,\\
\arg\left( v(h = 2) \right) &= \frac{4\pi}{- 3} = \frac{2\pi}{3} = 120^\circ,\\
\arg\left( v(k = 0) \right) &= \frac{\pi}{- 2} =  270^\circ,\\
\arg\left( v(k = 1) \right) &= \frac{3\pi}{- 2} = \frac{\pi}{2} = 90^\circ
\end{align*}

Only $\arg\left( v(h = 2) \right) = 120^\circ$ falls within the sector $[108^\circ, 180^\circ]$.

\textbf{Key observation}: The original series for branch $n=2$:
\[
y_{2} = vM\left( 1;5;2;7v^{2} \right),\quad v = e^{\frac{4\pi i}{5}}
\]
\textbf{diverges} for $x=7$ since $|x| = 7 > R$, where $R$ is the convergence radius.

\textbf{Conclusion}: The analytical continuation of the \textbf{divergent series} for branch $y_2$ is the \textbf{convergent series}:
\[
y(h = 2) = vM\left( 1; - 3; - 5;\frac{v^{- 5}}{2} \right),\quad v = e^{\frac{\ln\left| \frac{1}{35} \right| + 4\pi i}{- 3}}
\]
where
\begin{equation}
M(m; \alpha; \beta; z) = m + z + \sum_{\ell=2}^{\infty} \frac{z^{\ell}}{\ell!} \prod_{\gamma=1}^{\ell-1} (m - \alpha \gamma + \beta \ell)
\end{equation}
This continuation yields the same root as the original branch $y_2$ but through a series that converges for $|x| \geq R$.
The convergence radius for this case is:
\[
R = \frac{|1 - 5/2|^{\frac{2}{5}}}{|2 - 5|} = \frac{|1 - 2.5|^{0.4}}{3} \approx \frac{1.5^{0.4}}{3} \approx 0.47
\]
Since $|x| = 7 > 0.47$, the original series diverges.

\begin{table}[ht]
\centering
\caption{Branch trajectories of the ultra-radical $\sqrt[{n;5;2}]{x}$ as $x$ varies from -2 to 2}
\begin{tabular}{cccc}
\hline
\textbf{$n=1$} & \textbf{$n=2$} & \textbf{$n=-2$} & \textbf{$n=-1$} \\
\includegraphics[width=0.22\textwidth]{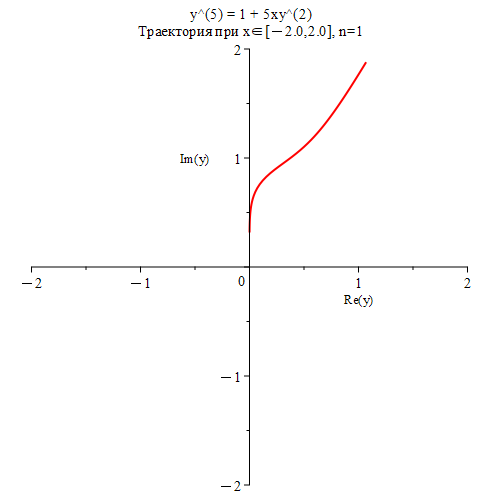} & 
\includegraphics[width=0.22\textwidth]{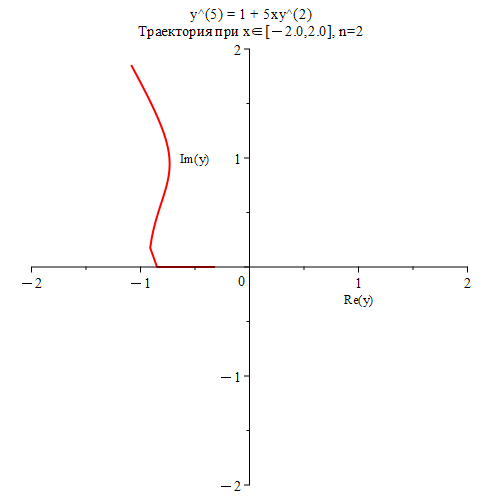} &
\includegraphics[width=0.22\textwidth]{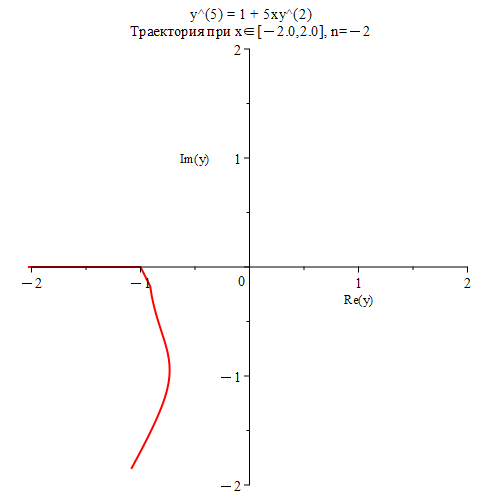} &
\includegraphics[width=0.22\textwidth]{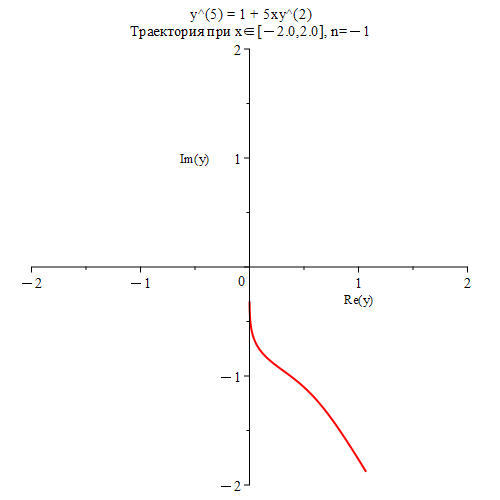} \\ 
\hline
\end{tabular}
\end{table}

\textbf{Remarks:}
\begin{itemize}
    \item The branch $n=0$ of the ultra-radical $\sqrt[{0;5;2}]{x}$ for real $x$ always returns a positive real $y$; therefore its plot on the complex plane (a straight line) is not shown here.
    \item Unlike the ordinary fifth root $\sqrt[5]{x}$, the branches $n=2$ and $n=-2$ of the ultra-radical possess segments with negative real values (visible in the plots).
    \item For $|x|<R$ (inside the convergence radius) all branches are given by the original master series. For $|x|\geq R$ analytic continuation via conjugate series is required, as demonstrated in the solution above for $n=2$.
\end{itemize}

\subsection{Example 2: Candidates L(h) and L(k) on Sector Boundaries}
Consider the ultra-radical $\sqrt[{n;4;1}]{x}$, solving the equation $y^{4} = 1 + 4xy$ for parameter values where two roots approach their intersection point. Let $x = P R \exp\left(Q\frac{\pi i}{4}\right)$ with varying parameters $P$ and $Q$.
\begin{table}[ht]
\centering
\caption{Root behavior for different parameter combinations}
\label{tab:root_behavior}
\begin{tabular}{ccc}
\hline
$P = 0.9,\ Q = 0.9$ & $P = 0.9,\ Q = 1$ & $P = 0.9,\ Q = 1.1$ \\
\includegraphics[width=0.3\textwidth]{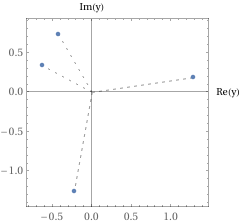} &
\includegraphics[width=0.3\textwidth]{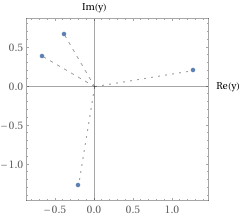} &
\includegraphics[width=0.3\textwidth]{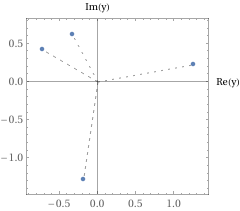} \\
\hline
$P = 1.0,\ Q = 0.9$ & $P = 1.0,\ Q = 1$ & $P = 1.0,\ Q = 1.1$ \\
\includegraphics[width=0.3\textwidth]{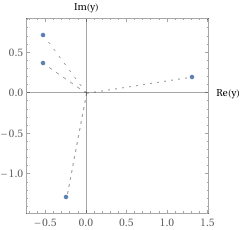} &
\includegraphics[width=0.3\textwidth]{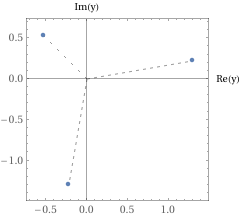} &
\includegraphics[width=0.3\textwidth]{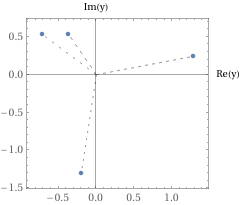} \\
\hline
$P = 1.1,\ Q = 0.9$ & $P = 1.1,\ Q = 1$ & $P = 1.1,\ Q = 1.1$ \\
\includegraphics[width=0.3\textwidth]{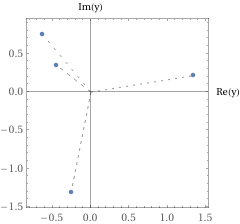} &
\includegraphics[width=0.3\textwidth]{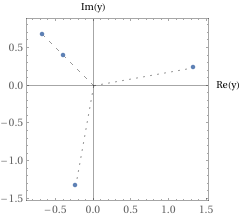} &
\includegraphics[width=0.3\textwidth]{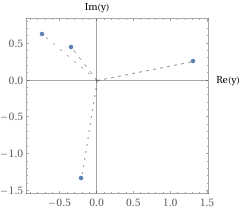} \\
\hline
\end{tabular}
\end{table}

\textbf{Observations}:
\begin{itemize}
\item \textbf{For $P < 1$}: All roots stably determined by original series within convergence radius
\item \textbf{For $P = 1$ and $Q = 1$}: Bifurcation point observed - two roots coincide
\item \textbf{For $P \geq 1$}: Branch redistribution occurs across sectors (see Example 1)
\item \textbf{For $P > 1$ and $Q = 1$}: In this degenerate case at sector boundaries, where multiple candidates for continuation exist, we propose the following convention: when two branches $h$ and $k$ intersect on the boundary of sectors $n-1$ and $n$, the branch with index $n$ is continued via the $h$-series, and the branch with index $n-1$ via the $k$-series. A rigorous proof of this convention and its general statement are given in~\cite{Berezin2026_ULog} (Theorem~1). This ensures synchronization of the numbering of branches $0$ and $-1$ of the ultraradical with branches $W_0$ and $W_{-1}$ of the Lambert $W$-function.

\end{itemize}

The graphics show root trajectories on the complex plane, demonstrating:
\begin{itemize}
\item Continuous evolution of roots as parameters vary
\item The geometric criterion successfully maintains branch continuity
\end{itemize}

\section{Visualization of Ultra-Radicals}

The ultra-radicals $\sqrt[n;a;b]{x}$ exhibit rich geometric structures across different parameter combinations. Table~\ref{tab:ultra_comparison} shows four distinct cases:

\begin{table}[ht]
\centering
\begin{tabular}{cccc}
\includegraphics[width=0.22\textwidth]{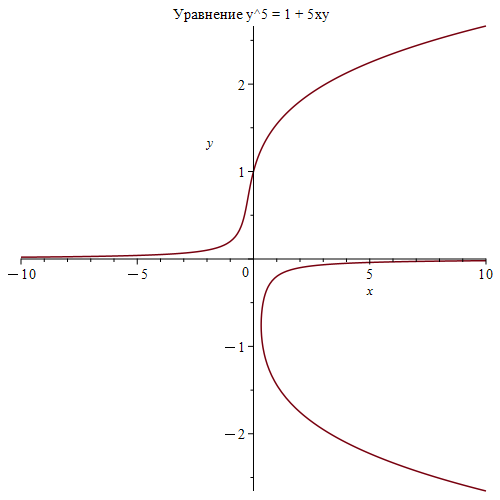} & 
\includegraphics[width=0.22\textwidth]{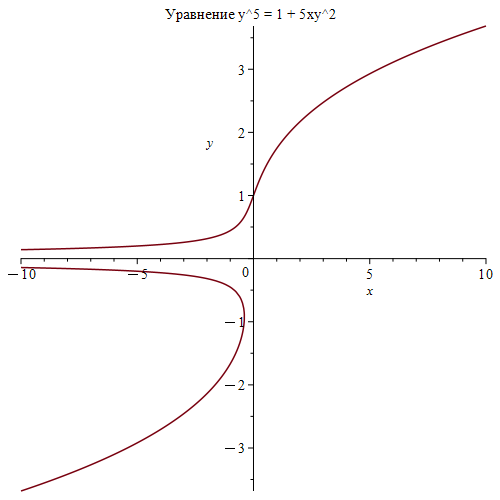} &
\includegraphics[width=0.22\textwidth]{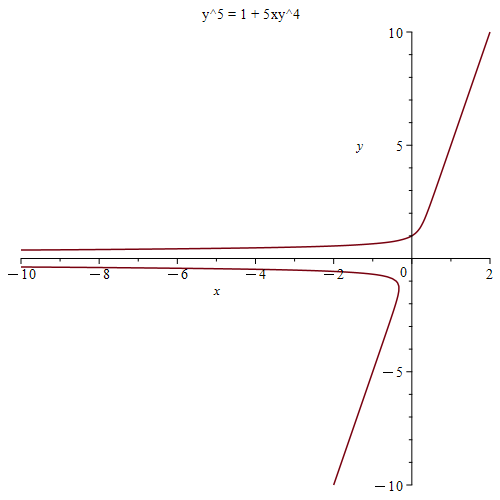} &
\includegraphics[width=0.22\textwidth]{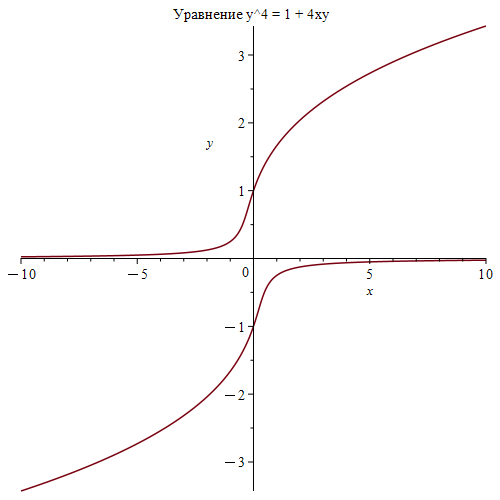} \\ 
\end{tabular}
\caption{Comparison of ultra-radicals with different parameters. Note the consistent behavior of the $n=0$ branch across all equations.}
\label{tab:ultra_comparison}
\end{table}

\subsection{Special Case: n=0 Branch}

For real-valued $x$, the branch $n=0$ of the ultra-radical $\sqrt[0;a;b]{x}$ exhibits distinctive properties that make it particularly valuable in applications:

\begin{itemize}
    \item \textbf{Real and Positive}: The $n=0$ branch always returns real positive values for real $x$.
    
    \item \textbf{Deterministic Series Selection}: For $|x| \geq R$ (outside convergence radius), the geometric criterion automatically selects:
    \begin{itemize}
        \item $h=0$ series for $x > 0$
        \item $k=0$ series for $x < 0$
    \end{itemize}
    For $|x| < R$, the original series with $n=0$ is used.
\end{itemize}

These properties make the $n=0$ branch the natural choice for most practical applications, similar to how the principal value is preferred for $\sqrt{x}$ and $\ln x$ in conventional analysis.

Maple code for generating these plots:
\begin{verbatim}
plots:-implicitplot(y^a=1+a*x*y^b,x=-10..10,y=-10..10,title="y^a = 1 + axy^b");
\end{verbatim}

\subsection{Structural Stability of the Principal Branch}

The special behavior of the $n=0$ branch reveals a profound structural property: 
among all branches, only this one evolves \emph{smoothly} as the parameters $a$ and $b$ vary continuously.

\noindent\textbf{Observation:} When $a$ is varied continuously from, say, 3 to 12, 
the $n=0$ branch changes gradually — its graph simply adjusts its curvature and asymptotic behavior. 
All other branches, however, appear and disappear abruptly as the number of solutions changes 
discontinuously with $a$. For irrational $a$, the equation technically possesses infinitely many 
branches, yet only the $n=0$ branch remains stable and well‑defined across the entire parameter range.

\noindent\textbf{Physical and computational implications:}
\begin{itemize}
    \item In physical systems where exponents $a$, $b$ themselves vary with time or state 
    (e.g., non‑linear media with field‑dependent exponents, adaptive systems), 
    only the principal branch can be tracked continuously.
    
    \item Iterative methods (Newton, Halley, etc.) have no inherent mechanism to identify 
    which of the infinitely many numerical roots corresponds to the principal branch — 
    they simply converge to whichever root lies closest to the initial guess.
    
    \item The Master‑J method, in contrast, \emph{identifies} the principal branch explicitly: 
    it is the one for which the phase factor $e^{2\pi i n / a}$ reduces to unity ($n=0$), 
    eliminating the oscillatory term that complicates all other branches.
\end{itemize}

Thus, for the principal branch, the solution simplifies to:
\[
\sqrt[0;a;b]{x} = M\!\left(1; a; b; x\right) \quad (\text{for }|x|<R),
\]
with the analytic continuation automatically selecting $h=0$ (for $x>0$) or $k=0$ (for $x<0$) when $|x|\ge R$. No other branch enjoys such algorithmic simplicity.

In problems where the argument exceeds the radius of convergence of the original series (e.g., in thermodynamic and plasma physics equations of the form $\Phi = \mathcal{K}(1-w^{-2}) + \Theta(1-w^{\gamma-1})$), stationary solutions are governed by two primary branches of analytic continuation, denoted as $h=0$ and $k=0$. Within a unified indexing scheme consistent for both the ultraradical and the Lambert $W$ function, these branches correspond to indices $n=0$ and $n=-1$, respectively. As a rule, the physically admissible solution belongs to the branch satisfying $w \to 1$ as $\Phi \to 0$, which corresponds to the system's return to its unperturbed equilibrium state.

\subsection{Geometric Criterion for the Principal Branch (n=0)}

For the principal branch $y_0(x)=\sqrt[0;a;b]{x}$, the branch‑selection algorithm simplifies drastically. When $|x|\ge R$ (outside the convergence radius), the branch $n=0$ always admits two continuation candidates: $h=0$ and $k=0$. The choice between them is governed by a single geometric inequality:

\[
\boxed{%
\begin{aligned}
&\text{Use the $h$-series with } h=0 \quad \text{if } |\arg(x)| \le \frac{\pi\,|b-a|}{a},\\[6pt]
&\text{otherwise use the $k$-series with } k=0.
\end{aligned}}
\]

\section{Generalizations}
\subsection{Equations with Arbitrary Coefficients}
For equations with arbitrary coefficients:
\[
pu^{a} = q + zu^{b}
\]
we apply the substitution:
\[
u = y\left( \frac{q}{p} \right)^{\frac{1}{a}},\quad z = aqx\left( \frac{p}{q} \right)^{\frac{b}{a}}
\]
to obtain the canonical form:
\[
y^{a} = 1 + axy^{b}
\]

\subsection{Equations with Arbitrary Number of Terms}

For equations with multiple terms:
\[
py^{a} = q + x_{1}y^{b_{1}} + x_{2}y^{b_{2}} + \cdots
\]
such equations are solved by merging several ultra-radicals:
\[
y = v \cdot _{a}^{1}:_{\frac{x_{1}v^{b_{1}}}{aq},\ \frac{x_{2}v^{b_{2}}}{aq},\ \cdots}^{b_{1},\ b_{2},\ \cdots},\quad v = \sqrt[a]{\left| \frac{q}{p} \right|}e^{\frac{\left( \arg\left( \frac{q}{p} \right) + 2\pi n \right)i}{a}},\ n \in \mathbb{Z}
\]

The operation of merging master series will be detailed in the appendices.

\section{The Ultra‑Radical as a Solution of Nonlinear ODEs}
\label{sec:ultra_ode}
The differential properties of the ultra‑radical reveal its important role in the theory of nonlinear ordinary differential equations.

Let $y = \sqrt[{n;a;b}]{x}$ be the ultra‑radical defined by the equation $y^{a} = 1 + a x y^{b}$. Differentiating this identity with respect to $x$ yields
\[
a y^{a-1} y' = a y^{b} + a b x y^{b-1} y'.
\]
Cancelling $a$ (assuming $a \neq 0$) and solving for $y'$, we obtain the compact algebraic formula

\begin{equation}
\boxed{\,
\frac{d}{dx}\sqrt[{n;a;b}]{x} 
= \frac{(\sqrt[{n;a;b}]{x})^{b-a+1}}{1 - b x (\sqrt[{n;a;b}]{x})^{b-a}}
= \frac{\sqrt[{n;a;b}]{x}}{(\sqrt[{n;a;b}]{x})^{a-b} - b x}.
\,}
\label{eq:ultra_derivative}
\end{equation}

Thus the ultra‑radical satisfies the nonlinear differential equation

\begin{equation}
\boxed{\,
y' = \frac{y^{b-a+1}}{1 - b x y^{b-a}},
\qquad y(0)=1.
\,}
\label{eq:ultra_ode_en}
\end{equation}

Conversely, any ODE that can be reduced to the form \eqref{eq:ultra_ode_en} has the solution $y(x) = \sqrt[{n;a;b}]{x}$.

\subsection{Practical Use}
\begin{itemize}
    \item \textbf{Computing the derivative:} Once the value of the ultra‑radical is known, its derivative follows directly from \eqref{eq:ultra_derivative}.
    \item \textbf{Solving ODEs:} Equations of the form \eqref{eq:ultra_ode_en} are solved explicitly via the ultra‑radical.
    \item \textbf{Verification of numerical methods:} The ultra‑radical provides a benchmark analytical solution for testing numerical ODE solvers.
    \item \textbf{Behaviour analysis:} Formula \eqref{eq:ultra_derivative} allows one to study monotonicity, convexity and asymptotic behaviour of the solution.
\end{itemize}

\section{Integral of the Ultra-Radical}
\label{sec:ultra_integral}

Similar to the derivative, the indefinite integral of the ultra-radical can also be expressed through the ultra-radical itself. This extends the analytic capabilities of the method and enables, in particular, the derivation of explicit expressions for pseudopotentials in nonlinear wave problems. The integral is most conveniently obtained by the method of integration of the inverse function, since the inverse relation \(x = x(y)\) for the ultra-radical has a simple power‑law form.

\subsection{Method of integration via the inverse function}

If \(y = y(x)\) is an invertible function and \(x = x(y)\) is its inverse, the identity
\[
\int y(x)\,dx = x y(x) - \int x(y)\,dy
\]
holds; it follows from integration by parts or directly from geometric considerations (see, e.g., standard textbooks of mathematical analysis). This method is particularly efficient when \(x(y)\) has a simpler expression than the original dependence \(y(x)\).

\subsection{Derivation of the integral formula}

Let \(y = \sqrt[{n;a;b}]{x}\) be the ultra‑radical satisfying the defining equation
\[
y^{a} = 1 + a x y^{b}, \qquad a \neq 0,\; b \neq 1.
\]
Solve it for \(x\) as a function of \(y\):
\begin{equation}
x = \frac{y^{a-b} - y^{-b}}{a}. \label{eq:x_from_y_integral_en}
\end{equation}
Now apply the integration‑via‑inverse‑function formula:
\[
\int y\,dx = x y - \int x\,dy.
\]
Substituting \(x\) from \eqref{eq:x_from_y_integral_en} we obtain
\begin{align*}
\int y\,dx &= x y - \frac{1}{a}\int \bigl( y^{a-b} - y^{-b} \bigr)\,dy \\
&= x y - \frac{1}{a}\Biggl[ \frac{y^{a-b+1}}{a-b+1} - \frac{y^{1-b}}{1-b} \Biggr] + C,
\end{align*}
where \(C\) is an arbitrary constant of integration, and we assume \(a-b+1 \neq 0\) and \(b \neq 1\) (the case \(b=1\) is treated separately).

Thus we arrive at the general formula for the indefinite integral:
\begin{equation}
\boxed{\,
\int \sqrt[{n;a;b}]{x}\,dx 
= x\sqrt[{n;a;b}]{x} 
- \frac{1}{a}\Biggl[ 
\frac{(\sqrt[{n;a;b}]{x})^{a-b+1}}{a-b+1} 
- \frac{(\sqrt[{n;a;b}]{x})^{1-b}}{1-b} 
\Biggr] + C.
\label{eq:ultra_integral_general}
\,}
\end{equation}

\subsection{Normalized form used in physical applications}

In many physical problems one requires the integral to vanish at \(x=0\), which corresponds to the condition \(y(0)=1\). Substituting \(x=0\), \(y=1\) into \eqref{eq:ultra_integral_general} fixes the constant \(C\) and yields a formula for the definite integral from \(0\) to \(x\):
\begin{equation}
\boxed{\,
\int_{0}^{x} \sqrt[{0;a;b}]{t}\,dt 
= \frac{1}{a}\Biggl[ 
\frac{a-b}{a-b+1}\Bigl( y^{a-b+1} - 1 \Bigr) 
+ \frac{b}{1-b}\Bigl( y^{1-b} - 1 \Bigr) 
\Biggr], 
\qquad y = \sqrt[{0;a;b}]{x}.
\label{eq:ultra_integral_zero_en}
\,}
\end{equation}
Exactly this normalized form arises naturally when deriving a pseudopotential from Poisson’s equation after substituting the concentration expressed via the ultra‑radical.

\subsection{Special case b=1}

For \(b=1\) formula \eqref{eq:ultra_integral_general} becomes
\begin{equation}
\int \sqrt[{n;a;1}]{x}\,dx 
= x\sqrt[{n;a;1}]{x} 
- \frac{1}{a}\Biggl[ 
\frac{(\sqrt[{n;a;1}]{x})^{a}}{a} 
- \ln\!\bigl( \sqrt[{n;a;1}]{x} \bigr) 
\Biggr] + C,
\label{eq:ultra_integral_b1_en}
\end{equation}
which can be obtained either by direct computation or by taking the limit \(b \to 1\) in \eqref{eq:ultra_integral_general}.

\subsection{Verification and an example}

To verify, differentiate the right‑hand side of \eqref{eq:ultra_integral_general} using the derivative formula for the ultra‑radical \eqref{eq:ultra_derivative} and the defining equation. After algebraic simplification one confirms that the derivative equals \(\sqrt[{n;a;b}]{x}\), which validates the integral.

As an illustration take \(a=3\), \(b=2\) (the equation \(y^{3}=1+3xy^{2}\)). From \eqref{eq:ultra_integral_general}:
\[
\int \sqrt[{0;3;2}]{x}\,dx 
= x y - \frac{y^{2}}{6} + \frac{y^{-1}}{3} + C, \qquad y = \sqrt[{0;3;2}]{x}.
\]

\subsection{Conclusion}

The obtained formulas allow one to express integrals of solutions to a wide class of power‑law equations in closed analytic form. This opens the possibility of direct symbolic analysis in problems where previously a parametric representation or numerical integration was required.

\section{General principle: from conservation law to trinomial}
\label{sec:general_principle}

In numerous stationary physical processes, the key principle is the conservation of the total energy per elementary volume or particle:
\begin{equation}
E_{\text{kin}} + E_{\text{int}} + E_{\text{pot}} = \text{const}.
\label{eq:energy_const}
\end{equation}
If we consider two states of the system (initial state $1$ and current state $2$) and subtract their energy balances, the constant cancels out, and we obtain the zero-sum condition for the increments:
\begin{equation}
\Delta E_{\text{kin}} + \Delta E_{\text{int}} + \Delta E_{\text{pot}} = 0.
\label{eq:energy_delta}
\end{equation}
For a broad class of media, kinetic and internal energies depend on characteristic state variables (velocity, density, temperature) via power-law relations. After subtracting the initial state, each increment takes the form $A(1 - \alpha)$, where $\alpha$ is the ratio of the current value to the initial one. Introducing a universal dimensionless variable $w$ that relates these ratios through the continuity equation or equation of state, Eq.~\eqref{eq:energy_delta} reduces to the canonical form:
\begin{equation}
\mathcal{K}\left(1 - \frac{1}{w^2}\right) + \frac{\Theta}{\gamma-1}\left(1 - w^{\gamma-1}\right) - \Phi = 0.
\label{eq:general_trinomial}
\end{equation}
Here $\mathcal{K}$, $\Theta$, and $\Phi$ are energy scales\footnote{In master series formulas, $\gamma$ denotes the summation index; in physical applications, the adiabatic index. The context uniquely determines the meaning.} determined by the physics of the process, and $w$ is a universal unknown describing the degree of deviation of the system from equilibrium. In the applications considered, $w$ represents a ratio of physical quantities (densities, velocities, concentrations), so physically $w > 0$. Multiplying Eq.~\eqref{eq:general_trinomial} by $w^2$ immediately yields a trinomial equation solvable via the ultraradical. For such equations, the ultraradical yields positive real solutions on branches $n=-1$ and $n=0$. The physically admissible solution is identified with the branch satisfying $w \to 1$ as $\Phi \to 0$, which corresponds to the system's transition to its unperturbed equilibrium state. Below, we demonstrate how this general algorithm is implemented in gas dynamics and plasma physics.

\subsection{Example 1: Gas dynamics in a gravitational field (derivation sketch)}
\label{subsec:gas_derivation}

Consider a stationary flow of an ideal gas in a gravitational field within a pipe of constant cross-section. The energy conservation law for an elementary mass layer $m_s$ reads:
\begin{equation}
\frac{m_s u^2}{2} + \frac{\gamma}{\gamma-1} N_s k_B T + m_s g h = \text{const}.
\end{equation}
Subtracting the inlet state ($u_1, T_1, h_1$) and using the continuity equation 
$\rho_1 u_1 = \rho_2 u_2 \Rightarrow u_2 = u_1 w^{-1}$, where $w = \rho_2/\rho_1 = u_1/u_2$, 
together with the adiabatic relation $T_2 = T_1 w^{\gamma-1}$, we obtain:
\begin{equation}
\frac{m_s u_1^2}{2}\left(1 - \frac{1}{w^2}\right) + 
\frac{\gamma k_B T_1 N_s}{\gamma-1}\left(1 - w^{\gamma-1}\right) + 
m_s g (h_1 - h_2) = 0.
\end{equation}
Introducing the energy variables $\mathcal{K} = \frac{1}{2}m_s u_1^2$, 
$\Theta = \gamma k_B T_1 N_s$, $\Phi = m_s g H$, 
we arrive exactly at the form~\eqref{eq:general_trinomial}. The critical parameter $H_{\max}$ is found from the condition $|x|=R$ for the ultraradical, which physically corresponds to reaching the local speed of sound $u_2 = c_s$.

\subsection{Example 2: Plasma physics (brief derivation via isomorphism)}
\label{subsec:plasma_derivation}

In a stationary ion-acoustic wave in plasma, the ion energy balance in the co-moving frame $\xi = z - Vt$ has the same structure~\eqref{eq:general_trinomial}, but with replaced physical scales:
\begin{equation}
\mathcal{K} = \frac{m_{i0} V^2}{2}, \quad 
\Theta = \gamma_i k_B T_{i0}, \quad 
\Phi = Z e \varphi.
\end{equation}
The continuity equation yields $v = -V/\eta$, where $\eta = n/n_0$ is the density ratio. The adiabatic ion state relates $T_i = T_{i0} \eta^{\gamma_i-1}$. Substitution into~\eqref{eq:general_trinomial} gives:
\begin{equation}
\frac{m_{i0} V^2}{2}\left(1 - \frac{1}{\eta^2}\right) + 
\frac{\gamma_i k_B T_{i0}}{\gamma_i-1}\left(1 - \eta^{\gamma_i-1}\right) - 
Z e \varphi = 0.
\label{eq:plasma_final}
\end{equation}
The mathematical structure of~\eqref{eq:plasma_final} is identical to the gas-dynamic case: the gravitational potential $m_s g H$ is replaced by the electrostatic potential $Z e \varphi$, and the density ratio $\rho_2/\rho_1$ by the concentration ratio $\eta = n/n_0$. Consequently, all analytical results obtained for the gas (solution via ultraradical, radius of convergence, critical condition $d\Phi/d\eta=0$) automatically transfer to the plasma case.

\subsection{Nonlinear electrical circuits and variable exponent}
\label{sec:other_applications}

Although the ultraradical $\sqrt[{n;a;b}]{x}$ was introduced as a solution to the algebraic-transcendental equation $y^a = 1 + a x y^b$, its mathematical structure turns out to be naturally connected with the description of a broad class of physical systems --- in particular, electrical circuits with nonlinear elements. One of the most illustrative examples is the analysis of nonlinear electrical circuits, where voltage-current characteristics of components are often described by power-law functions with fractional or even irrational exponents. In this subsection we show how the ultraradical naturally arises in circuit theory and how it allows correct treatment of even the infinite multivaluedness of solutions.

\subsubsection{Power-law nonlinearity and its physical meaning}

For a linear resistor, the relation between voltage $U$ and current $I$ is given by Ohm's law:
\[
U = R I.
\]
However, for many real devices (incandescent lamps, semiconductors, varistors) this dependence is essentially nonlinear. It is often approximated by a power-law function
\begin{equation}
U = \alpha I^{k}, \qquad k > 0,\; k \neq 1,
\label{eq:power_law}
\end{equation}
where $\alpha > 0$ is a coefficient and the exponent $k$ characterizes the degree of nonlinearity. For $k>1$ the resistance increases with current (incandescent lamp), for $k<1$ it decreases (thermistor). In the case $k=1$ one recovers the ordinary linear resistor.

Consider the simplest circuit: a voltage source $U_0$, a linear resistor $R_0$, and a nonlinear element with characteristic~\eqref{eq:power_law}. By Kirchhoff's voltage law:
\[
U_0 = R_0 I + \alpha I^{k}.
\]
After moving all terms to one side we obtain a trinomial equation:
\begin{equation}
\alpha I^{k} + R_0 I - U_0 = 0.
\label{eq:circuit_eq}
\end{equation}

\subsubsection{Problem of fractional exponent and multivaluedness}

If the exponent $k$ is an irrational number, Eq.~\eqref{eq:circuit_eq} possesses an \emph{infinite} set of complex solutions. This is a direct consequence of the multivaluedness of the power function $I^{k}$ for non-integer $k$.

The ultraradical $\sqrt[{n;a;b}]{x}$ explicitly parametrizes all solutions by an integer index $n \in \mathbb{Z}$, which enables:
\begin{enumerate}
    \item identification of each solution as belonging to a specific branch;
    \item tracking of continuous solution variation under parameter changes (including the exponent $k$ itself);
    \item performing analytical operations (differentiation, integration) with full branch control.
\end{enumerate}

\subsubsection{Variable exponent: thermal drift}

In real devices the exponent $k$ may vary with time, e.g., due to heating. Let $k = k(t)$ be a slowly varying function. Then Eq.~\eqref{eq:circuit_eq} becomes parametric, and its solution becomes a function of time:
\[
I(t) = \left( \frac{U_0}{\alpha} \right)^{1/k(t)} \cdot \sqrt[{0;\,k(t);\,1}]{x(t)},
\]
where
\[
x(t) = - \frac{R_0}{k(t)\,U_0} \left( \frac{U_0}{\alpha} \right)^{1/k(t)}.
\]

Thanks to the continuity of the principal branch of the ultraradical, the function $I(t)$ remains smooth under any smooth variations of $k(t)$.

\section{The Ultralogarithm: The Logarithmic Core of the Theory}

The presented theory of the ultra radical possesses a profound and symmetric "logarithmic" counterpart. This counterpart not only simplifies the formalism but also directly bridges it with classical analysis, revealing it as a natural generalization of familiar functions.

\subsection{From Ultra Radical to Ultralogarithm}

The fundamental identity linking master series, $M(1;a;b;x) = \exp(M(0;a;b;x))$, suggests that $M(0;a;b;x)$ is the "logarithm" of the ultra radical. Crucially, this is not merely a formal composition---it is an \textbf{independent canonical function} solving its own defining equation. For $a=1$, this equation takes an elegantly simple form, which we term the \textbf{$b$-logarithm equation}:
\begin{equation}
\label{eq:b_log_def_en}
\boxed{\, y = \ln\!\bigl(1 + x e^{\,b y}\bigr) \,}
\end{equation}
Its solution, $y = \operatorname{ulog}_b(x) := M(0;1;b;x)$, is called the \textbf{$b$-logarithm} (or ultralogarithm).

This reveals \emph{two parallel families} of canonical functions:
\begin{center}
\renewcommand{\arraystretch}{1.2}
\begin{tabular}{c|c|c}
\textbf{Type} & \textbf{General Equation} & \textbf{Canonical Form ($a=1$)} \\ \hline
\textbf{Power-type} & $y^a = 1 + a x y^b$ & $y = (1 + x y^b)^{1}$ (trivial) \\
\textbf{Logarithmic-type} & $y = \dfrac{\ln(1 + a x e^{b y})}{a}$ & $\boxed{y = \ln(1 + x e^{b y})}$ \\
\end{tabular}
\end{center}

The key insight is that the **canonical logarithmic form is non-degenerate**. It defines a nontrivial one-parameter family $\operatorname{ulog}_b(x)$, which generalizes the natural logarithm: $\ln(1+x) = \operatorname{ulog}_0(x)$. Furthermore, specific parameter choices recover other classic functions. For instance, with $a=2, b=1$, the general logarithmic master function $M(0;2;1;x)$ yields the \textbf{inverse hyperbolic sine}: $\operatorname{arsinh}(x) = M(0;2;1;x)$. This is not a coincidence but a consequence of the underlying structure: the condition $a - b = 1$ (as seen in the graphs of the corresponding ultra radical) and $b=1$ (which simplifies integration formulas) mark this function as a distinguished, well-behaved node within the parametric space.

\subsection{Branch Relationships: From Ultra Radical to Ultralogarithm}

The defining equations of the ultra radical and the ultralogarithm form a coupled system via the logarithm operation:
\begin{align}
    y &= \bigl(1 + a x y^{b}\bigr)^{1/a} \label{eq:ultraradical_eq_en} \quad \text{(Ultra radical $\sqrt[{n;a;b}]{x}$)},\\
    u &= \frac{\ln\!\bigl(1 + a x e^{b u}\bigr)}{a} \label{eq:ultralogarithm_eq_en} \quad \text{(Ultralogarithm $M(0;a;b;x)$)}.
\end{align}
The substitution $y = e^{u}$ (or $u = \ln y$) transforms \eqref{eq:ultraradical_eq_en} into \eqref{eq:ultralogarithm_eq_en}. However, due to the multivalued nature of the complex logarithm, this connection generates the complete family of branches.

For the canonical form, the general solution is given by parametric formulas:
\begin{align}
    y_n &= v_n \cdot M\!\bigl(1;\; a;\; b;\; x V_n\bigr), \label{eq:ultraradical_solution_en} \\
    u_n &= \ln(v_n) + M\!\bigl(0;\; a;\; b;\; x V_n\bigr), \label{eq:ultralogarithm_solution_en}
\end{align}
where
\[
v_n = \exp\!\left(\frac{\ln|1| + i\bigl(\arg(1) + 2\pi n\bigr)}{a}\right) = e^{2\pi i n / a}, \quad
V_n = v_n^{b} = e^{2\pi i b n / a}, \quad n \in \mathbb{Z}.
\]

\subsubsection*{Interpretation and Implications}
\begin{itemize}
    \item Formulas \eqref{eq:ultraradical_solution_en} and \eqref{eq:ultralogarithm_solution_en} demonstrate that the branch $u_n$ of the ultralogarithm is the exact logarithm of the corresponding branch of the ultraradical ($u_n = \ln y_n$). However, this relation is nontrivial: it is expressed not as a simple master series $M(0; a; b; x)$, but as a shift by a complex constant $\ln(v_n)$ and a scaling of the argument to $x V_n$ within the logarithmic master series.

    \item The expression $M(0;a;b;x V_n)$ in \eqref{eq:ultralogarithm_solution_en} represents the \textbf{power series for any branch of the ultralogarithm}. Thus, taking the logarithm of the ultra radical yields not an arbitrary transcendental function, but a specific, computable power series shifted by an integer imaginary constant.

    \item For $n=0$, we obtain the \emph{principal branches}: $v_0=1$, $V_0=1$, and the relationship simplifies to the expected one: $y_0 = M(1;a;b;x)$ and $u_0 = M(0;a;b;x)$, with $y_0 = e^{u_0}$.

    \item This formalism makes the \textbf{verification rule} explicit: for the substitution $y_n$ into equation \eqref{eq:ultraradical_eq_en} to yield an identity, it is necessary to compute $y_n^{a}$ as $\exp\bigl(a \cdot u_n\bigr)$, where $u_n$ is given by formula \eqref{eq:ultralogarithm_solution_en}. The integer $n$ in the expression for $v_n$ corresponds precisely to the verification parameter \texttt{u} computed automatically by the \texttt{SolverABC} algorithm.
\end{itemize}

This representation not only clarifies the multivalued structure but also provides a unified computational basis: an algorithm constructed for the ultra radical simultaneously computes the corresponding branches of the ultralogarithm via formulas \eqref{eq:ultraradical_solution_en}--\eqref{eq:ultralogarithm_solution_en}.

\subsection{Differential Properties and Relation to the Ultra Radical}
As follows from the exponential identity $M(1;a;b;x) = \exp(M(0;a;b;x))$, the derivatives of the master functions are linked by a fundamental relation:
\begin{equation}
\label{eq:master_diff_identity_en}
\frac{d}{dx}M(1;a;b;x) = M(1;a;b;x) \cdot \frac{d}{dx}M(0;a;b;x).
\end{equation}
This relationship reflects the logarithmic nature of the connection between the power-type and logarithmic-type series. For the $b$-logarithm ($a=1$), equation \eqref{eq:master_diff_identity_en} takes a concrete and illustrative form. Introducing $u = e^{\operatorname{ulog}_b(x)} = M(1;1;b;x)$, which is equivalent to $u = 1 + x u^b$, one can derive an explicit formula for the derivative that avoids cumbersome expressions:

\begin{equation}
\label{eq:b_log_derivative_explicit_en}
\frac{d}{dx} \operatorname{ulog}_b(x) = \frac{u^{\,b}}{b + (1 - b) u}.
\end{equation}

\textbf{Key observation:} Formula \eqref{eq:b_log_derivative_explicit_en}, unlike the direct expression for the derivative of the corresponding ultra radical, \emph{contains no explicit dependence on the independent variable $x$}. It expresses the rate of change of the $b$-logarithm solely in terms of its exponent $u$ and the parameter $b$, representing a significant algebraic simplification. This demonstrates that switching to the logarithmic parametrization does not merely alter the notation but \emph{simplifies the algebraic structure} of the fundamental relations.

\subsection{The Significance of the b-Logarithm}

\begin{enumerate}
    \item \textbf{Minimalism and Classical Bridge}. Like the ordinary logarithm, $\operatorname{ulog}_b(x)$ is defined by a single parameter $b$. It provides the most direct and intuitive entry point into the theory, seamlessly extending the familiar function $\ln(1+x)$.

    \item \textbf{Computational Advantages}. For small $|x|$, the series $M(0;1;b;x)$ often exhibits superior convergence and numerical stability compared to the power series of the corresponding ultra radical with parameters $(1, b)$.

    \item \textbf{A Transformational Tool}. Many equations involving exponentials can be reduced to the form solvable by the $b$-logarithm. Thus, $\operatorname{ulog}_b(x)$ acts as a \emph{universal solver} for a class of exponential-power equations, analogous to the role of the logarithm in linearizing multiplicative relationships.

    \item \textbf{A New Research Pathway}. The theory of the $b$-logarithm as a standalone object---its analytic properties, integrals, asymptotic behavior, and applications---constitutes a vast and nearly unexplored field. It offers a more accessible and focused domain for initial investigation than the full two-parameter theory of the ultra radical.
\end{enumerate}

\subsection{Methodological Impact}

Introducing the $b$-logarithm $\operatorname{ulog}_b(x)$ as a primary object serves a crucial purpose: it \textbf{grounds the theory in a simple, familiar core}. Instead of presenting the complete two-parameter apparatus immediately, one can first demonstrate how the well-known natural logarithm naturally extends by introducing a single parameter $b$. The subsequent generalization to the ultra radical (by introducing the parameter $a$) then appears as a logical and necessary step to capture the full algebraic structure, rather than an arbitrary abstraction.

Therefore, the ultralogarithm is not merely a corollary of the theory; it is its \textbf{conceptual keystone}. It facilitates understanding, strengthens connections to classical analysis, and provides a clear pathway for pedagogical exposition and practical application, all while highlighting profound connections to special functions like the inverse hyperbolic sine.

\section{Conclusion}

This work presents a systematic approach to solving the fundamental problem of analysis — ensuring continuity of branches of multi‑valued functions during analytical continuation beyond the radius of convergence of power series.

\textbf{Main results:}
\begin{enumerate}
    \item For the ultra-radical $\sqrt[{n;a;b}]{x}$ --- solution of the equation $y^{a} = 1 + a x y^{b}$ --- we propose a \textbf{geometric criterion} for branch selection. The criterion is based on the sector structure of the unit circle determined by parameter $a$, and guarantees root continuity when crossing the circle $|x| = R = \left| b^{-b/a} \cdot (b - a)^{(b-a)/a} \right|$.

    \item We developed and tested on examples a \textbf{deterministic algorithm} for analytical continuation, including standards for resolving ambiguities in degenerate cases (on sector boundaries).

    \item The framework \textbf{suggests a pathway for generalization} to equations with arbitrary coefficients and multiple terms through the merge operation, potentially extending its applicability to a broader class of equations.
\end{enumerate}

\appendix
\section{Derivation of Master Series}
\label{sec:master_series}

\subsection{Statement of the problem}

Consider the canonical unit equation:
\begin{equation}
y^a = 1 + axy^b,
\label{eq:canonical_1}
\end{equation}
where $a, b \in \mathbb{C}$ are complex exponents, $x \in \mathbb{C}$ is a complex variable.

We seek a solution of equation \eqref{eq:canonical_1} in the form of a power series:
\begin{equation}
y = 1 + \frac{N_1}{1!}x + \frac{N_2}{2!}x^2 + \frac{N_3}{3!}x^3 + \frac{N_4}{4!}x^4 + \cdots = 1 + \sum_{\ell=1}^{\infty} \frac{N_\ell}{\ell!} x^\ell,
\label{eq:series_ansatz}
\end{equation}
where $N_\ell$ are coefficients to be determined.

\subsection{Two-step expansion: generalized binomial and multinomial theorem}

To substitute series \eqref{eq:series_ansatz} into equation \eqref{eq:canonical_1}, we must expand the left-hand side $y^a$ as a power series in $x$. This process consists of two successive steps.

\textbf{Step 1: Generalized binomial theorem.}
Let $u = y - 1 = \sum_{\ell=1}^{\infty} \frac{N_\ell}{\ell!} x^\ell$. Then $y^a = (1 + u)^a$. 
We apply the generalized binomial theorem, valid for an arbitrary exponent $a \in \mathbb{C}$ for $|u| < 1$:
\begin{equation}
(1 + u)^a = \sum_{n=0}^{\infty} \binom{a}{n} u^n = 1 + \frac{a}{1!}u + \frac{a(a-1)}{2!}u^2 + \frac{a(a-1)(a-2)}{3!}u^3 + \frac{a(a-1)(a-2)(a-3)}{4!}u^4 + \cdots
\label{eq:binomial_expansion}
\end{equation}
where the generalized binomial coefficients are defined as
\begin{equation}
\binom{a}{n} = \frac{a(a-1)(a-2)\cdots(a-n+1)}{n!}.
\label{eq:binomial_coeff}
\end{equation}

\textbf{Step 2: Multinomial theorem.} 
Since $n$ in formula \eqref{eq:binomial_expansion} runs over non-negative integers $0, 1, 2, \ldots$, to expand each power $u^n$ we apply the classical multinomial theorem \cite{Graham1994}. We compute the first four powers:

\begin{align}
u^1 &= \frac{N_1}{1!}x + \frac{N_2}{2!}x^2 + \frac{N_3}{3!}x^3 + \frac{N_4}{4!}x^4 + \cdots, \label{eq:u1} \\
u^2 &= \left(\frac{N_1}{1!}x + \frac{N_2}{2!}x^2 + \frac{N_3}{3!}x^3 + \cdots\right)^2 = N_1^2 x^2 + N_1 N_2 x^3 + \left(\frac{N_2^2}{4} + \frac{N_1 N_3}{3}\right)x^4 + \cdots, \label{eq:u2} \\
u^3 &= \left(\frac{N_1}{1!}x + \frac{N_2}{2!}x^2 + \cdots\right)^3 = N_1^3 x^3 + \frac{3}{2}N_1^2 N_2 x^4 + \cdots, \label{eq:u3} \\
u^4 &= \left(\frac{N_1}{1!}x + \cdots\right)^4 = N_1^4 x^4 + \cdots. \label{eq:u4}
\end{align}

\subsection{Collecting coefficients}

Substituting expansions \eqref{eq:u1}--\eqref{eq:u4} into formula \eqref{eq:binomial_expansion} and grouping terms by powers of $x$, we obtain:

\begin{align}
(1+u)^a &= 1 + \frac{a}{1!}\left(\frac{N_1}{1!}x + \frac{N_2}{2!}x^2 + \frac{N_3}{3!}x^3 + \frac{N_4}{4!}x^4 + \cdots\right) \notag \\
&\quad + \frac{a(a-1)}{2!}\left(N_1^2 x^2 + N_1 N_2 x^3 + \left(\frac{N_2^2}{4} + \frac{N_1 N_3}{3}\right)x^4 + \cdots\right) \notag \\
&\quad + \frac{a(a-1)(a-2)}{3!}\left(N_1^3 x^3 + \frac{3}{2}N_1^2 N_2 x^4 + \cdots\right) \notag \\
&\quad + \frac{a(a-1)(a-2)(a-3)}{4!}\left(N_1^4 x^4 + \cdots\right) + \cdots
\end{align}

Collecting the coefficients for each power of $x$ yields:
\begin{align}
y^a &= 1 + aN_1 x + \left(\frac{aN_2}{2!} + \frac{a(a-1)N_1^2}{2!}\right)x^2 \notag \\
&\quad + \left(\frac{aN_3}{3!} + \frac{a(a-1)N_1 N_2}{2!} + \frac{a(a-1)(a-2)N_1^3}{3!}\right)x^3 \notag \\
&\quad + \left(\frac{aN_4}{4!} + \frac{a(a-1)}{2!}\left(\frac{N_2^2}{4} + \frac{N_1 N_3}{3}\right) + \frac{a(a-1)(a-2)N_1^2 N_2}{4} + \frac{a(a-1)(a-2)(a-3)N_1^4}{4!}\right)x^4 + \cdots
\label{eq:left_side}
\end{align}

\subsection{The right-hand side of the equation}

Now we expand the right-hand side of equation \eqref{eq:canonical_1}: $1 + axy^b$. 
Similarly to the left-hand side, we apply the generalized binomial theorem to $y^b = (1+u)^b$:
\begin{align}
y^b &= 1 + bN_1 x + \left(\frac{bN_2}{2!} + \frac{b(b-1)N_1^2}{2!}\right)x^2 \notag \\
&\quad + \left(\frac{bN_3}{3!} + \frac{b(b-1)N_1 N_2}{2!} + \frac{b(b-1)(b-2)N_1^3}{3!}\right)x^3 + \cdots
\label{eq:y_b}
\end{align}

Multiplying by $ax$, we get:
\begin{align}
1 + axy^b &= 1 + ax + abN_1 x^2 + a\left(\frac{bN_2}{2!} + \frac{b(b-1)N_1^2}{2!}\right)x^3 \notag \\
&\quad + a\left(\frac{bN_3}{3!} + \frac{b(b-1)N_1 N_2}{2!} + \frac{b(b-1)(b-2)N_1^3}{3!}\right)x^4 + \cdots
\label{eq:right_side}
\end{align}

\subsection{Equating coefficients}

Equating coefficients for equal powers of $x$ in expressions \eqref{eq:left_side} and \eqref{eq:right_side}, we obtain the system of equations:

\begin{align}
x^1: &\quad aN_1 = a, \label{eq:coeff1} \\
x^2: &\quad \frac{aN_2}{2!} + \frac{a(a-1)N_1^2}{2!} = abN_1, \label{eq:coeff2} \\
x^3: &\quad \frac{aN_3}{3!} + \frac{a(a-1)N_1 N_2}{2!} + \frac{a(a-1)(a-2)N_1^3}{3!} = a\left(\frac{bN_2}{2!} + \frac{b(b-1)N_1^2}{2!}\right), \label{eq:coeff3} \\
x^4: &\quad \frac{aN_4}{4!} + \frac{a(a-1)}{2!}\left(\frac{N_2^2}{4} + \frac{N_1 N_3}{3}\right) + \frac{a(a-1)(a-2)N_1^2 N_2}{4} + \frac{a(a-1)(a-2)(a-3)N_1^4}{4!} \notag \\
&\quad = a\left(\frac{bN_3}{3!} + \frac{b(b-1)N_1 N_2}{2!} + \frac{b(b-1)(b-2)N_1^3}{3!}\right). \label{eq:coeff4}
\end{align}

\subsection{Sequential solution of the system}

\textbf{Coefficient $N_1$.} 
From equation \eqref{eq:coeff1} we directly obtain:
\begin{equation}
N_1 = 1.
\label{eq:N1}
\end{equation}

\textbf{Coefficient $N_2$.} 
Substituting $N_1 = 1$ into equation \eqref{eq:coeff2} and dividing by $a/2!$, we get:
\begin{equation}
N_2 + (a-1) = 2b,
\end{equation}
hence:
\begin{equation}
\boxed{N_2 = 1 - a + 2b.}
\label{eq:N2}
\end{equation}

\textbf{Coefficient $N_3$.} 
Substituting $N_1 = 1$ into equation \eqref{eq:coeff3}, multiplying by $6/a$ and simplifying, we get:
\begin{equation}
N_3 + 3(a-1)N_2 + (a-1)(a-2) = 3bN_2 + 3b(b-1).
\end{equation}
Expressing $N_3$ and substituting $N_2 = 1-a+2b$:
\begin{align}
N_3 &= -3(a-1)(1-a+2b) - (a-1)(a-2) + 3b(1-a+2b) + 3b(b-1) \notag \\
&= -3a + 2a^2 - 9ab + 1 + 6b + 9b^2 \notag \\
&= (1-a+3b)(1-2a+3b).
\end{align}
Thus:
\begin{equation}
\boxed{N_3 = (1-a+3b)(1-2a+3b).}
\label{eq:N3}
\end{equation}

\textbf{Coefficient $N_4$.} 
Substituting $N_1 = 1$ into equation \eqref{eq:coeff4}, multiplying by $24/a$ and simplifying, we get:
\begin{align}
N_4 &= -12(a-1)\left(\frac{N_2^2}{4} + \frac{N_3}{3}\right) - 6(a-1)(a-2)N_2 - (a-1)(a-2)(a-3) \notag \\
&\quad + 4bN_3 + 12b(b-1)N_2 + 4b(b-1)(b-2).
\end{align}
Substituting expressions \eqref{eq:N2} and \eqref{eq:N3} for $N_2$ and $N_3$, after algebraic simplifications we obtain:
\begin{align}
N_4 &= -6a - 6a^3 - 96ab^2 + 11a^2 - 48ab + 44a^2b + 1 + 48b^2 + 12b + 64b^3 \notag \\
&= (1-a+4b)(1-2a+4b)(1-3a+4b).
\end{align}
Thus:
\begin{equation}
\boxed{N_4 = (1-a+4b)(1-2a+4b)(1-3a+4b).}
\label{eq:N4}
\end{equation}

\subsection{Master numbers and the master series}

We have obtained the coefficients $N_1, N_2, N_3, N_4$ of the single-valued power series \eqref{eq:series_ansatz}, which is a solution of the multivalued equation \eqref{eq:canonical_1}. 
Analyzing the structure of the obtained expressions, we observe that each coefficient $N_\ell$ is a product of linear factors of the form $(1 - a\gamma + b\ell)$, where $\gamma$ ranges from $1$ to $\ell-1$.

\textbf{Definition 1 (Master numbers).} 
We define the master number of order $\ell$ as the following product of factors (in the spirit of combinatorial notation \cite{Knuth1997}):
\begin{equation}
N_\ell = N(m; a; b; \ell) = \prod_{\gamma=1}^{\ell-1} (m - a\gamma + b\ell), \quad \ell \geq 2,
\label{eq:master_numbers_def}
\end{equation}
with $N_0 = 1$, $N_1 = 1$ (empty product).

For the canonical equation \eqref{eq:canonical_1}, the parameter $m = 1$, therefore:
\begin{align}
N_0 &= 1, \\
N_1 &= 1, \\
N_2 &= 1 - a + 2b, \\
N_3 &= (1-a+3b)(1-2a+3b), \\
N_4 &= (1-a+4b)(1-2a+4b)(1-3a+4b).
\end{align}

\textbf{Definition 2 (Master series).} 
We define the master series as the single-valued power series:
\begin{equation}
M(m; a; b; x) = 1 + \sum_{\ell=1}^{\infty} \frac{N(m; a; b; \ell)}{\ell!} x^\ell,
\label{eq:master_series_def}
\end{equation}
where $N(m; a; b; \ell)$ are the master numbers defined by formula \eqref{eq:master_numbers_def}.

For the canonical equation \eqref{eq:canonical_1}, the solution has the form:
\begin{equation}
y = M(1; a; b; x) = 1 + x + \frac{1-a+2b}{2!}x^2 + \frac{(1-a+3b)(1-2a+3b)}{3!}x^3 + \cdots
\label{eq:master_series_canonical}
\end{equation}

The parameter $m$ in the definition of master numbers is necessary for describing more general cases. In particular, there exists the identity $M(1; a; b; x) = \exp(M(0; a; b; x))$, relating master series with different values of the parameter $m$.

Furthermore, the product in formula \eqref{eq:master_numbers_def} can be expressed via the gamma function~\cite{Popov2017}:

\subsection{Expressing master numbers via the ratio of gamma functions}
The master number $N_\ell=\prod_{\gamma=1}^{\ell-1} \left( m - a \gamma + b \ell \right)$ can be expressed in closed form via the gamma function for $a \neq 0$ \cite{Popov2017}:

\begin{equation}
N_\ell = (-a)^{\ell-1} \, 
\frac{\Gamma\!\left( \dfrac{(a-b)\ell - m}{a} \right)}
     {\Gamma\!\left( \dfrac{a - m - \ell b}{a} \right)},
\label{eq:M_L_gamma}
\end{equation}

or in a numerically stable form:
\begin{equation}
N_\ell = a^{\ell-1} \, 
\frac{\Gamma\!\left( \dfrac{m + \ell b}{a} \right)}
     {\Gamma\!\left( \dfrac{m + \ell b}{a} - \ell + 1 \right)}.
\label{eq:M_L_gamma_clean}
\end{equation}

The gamma function has poles at points $z = 0, -1, -2, \dots$. Consequently, if the argument of the gamma function in the denominator takes a non-positive integer value, the denominator becomes infinite, and the entire master number vanishes (provided that the numerator does not simultaneously become infinite). This provides an explicit analytical criterion for identifying zero terms in the power series.

\textbf{Fundamental limitation of the gamma representation.}
Formulas \eqref{eq:M_L_gamma} and \eqref{eq:M_L_gamma_clean} serve as an \emph{auxiliary tool} for computing master numbers, but \emph{not as their replacement}. The gamma function cannot replace master numbers for three reasons.

\textbf{First}, for $a = 0$, master numbers are powers of $(m + b\ell)^{\ell-1}$, and in principle cannot be expressed as a ratio of gamma functions in the manner described above.

\textbf{Second}, for arbitrary parameter values $(m, a, b)$, master numbers do not reduce to the gamma function or to any other known generalizations of factorials (the Pochhammer symbol, multifactorials). The gamma function allows one to \emph{compute} a master number where possible, but it does not describe its \emph{structure}.

\textbf{Third}, there exist specific parameter values for which the argument of the gamma function simultaneously becomes a pole in both the numerator and the denominator of formulas \eqref{eq:M_L_gamma} and \eqref{eq:M_L_gamma_clean}. In this case, an $\infty/\infty$ indeterminate form arises, and the gamma representation loses its meaning, whereas the master number, computed according to definition, remains a perfectly well-defined finite number.

\begin{remark}[On the choice of gamma representation form]
For $m=1$ (ultra-radical), the numerically stable form \eqref{eq:M_L_gamma_clean} is preferred, as it avoids poles at typical parameter values. However, for $m=0$ (ultra-logarithm), form \eqref{eq:M_L_gamma_clean} may lead to an $\infty/\infty$ indeterminate form (for example, when $b=0$, both numerator and denominator simultaneously become poles). In this case, the historical form \eqref{eq:M_L_gamma} should be used, as it remains well-defined.
\end{remark}

\subsection{Connection with the confluent Fox--Wright function}
\label{subsec:fox_wright}

The confluent Fox--Wright function ${}_1\Psi_1$ is defined by two pairs of parameters and a single variable:
\begin{equation}
{}_1\Psi_1\!\left[\begin{matrix} (\alpha, A) \\ (\beta, B) \end{matrix}; z\right]
= \sum_{n=0}^{\infty} \frac{\Gamma(\alpha + An)}{\Gamma(\beta + Bn)}\, \frac{z^n}{n!}.
\label{eq:fox_wright_def}
\end{equation}

Substituting the gamma representation of master numbers \eqref{eq:M_L_gamma_clean} into the definition of the master series \eqref{eq:master_series_def}, we obtain:
\begin{equation}
M(m; a; b; x) = m + \sum_{\ell=1}^{\infty}
\frac{a^{\ell-1}}{\ell!}\,
\frac{\Gamma\!\left(\dfrac{m}{a} + \dfrac{b}{a}\,\ell\right)}
     {\Gamma\!\left(\dfrac{m}{a} + \left(\dfrac{b}{a}-1\right)\ell + 1\right)}\,
x^\ell.
\label{eq:master_gamma_series}
\end{equation}

Comparing \eqref{eq:master_gamma_series} with \eqref{eq:fox_wright_def}, we obtain the exact parameter correspondence:
\begin{equation}
\alpha = \frac{m}{a}, \qquad
A = \frac{b}{a}, \qquad
\beta = \frac{m}{a} + 1, \qquad
B = \frac{b}{a} - 1, \qquad
z = ax.
\label{eq:parameter_mapping}
\end{equation}

Since ${}_1\Psi_1$ starts at $n=0$, whereas the master series starts at $\ell=1$ (plus the constant term $m$), after elementary transformations we obtain the \textbf{identity}:
\begin{equation}
M(m; a; b; x) = m + \frac{1}{a}\left\{
{}_1\Psi_1\!\left[\begin{matrix} \left(\dfrac{m}{a},\; \dfrac{b}{a}\right) \\[6pt]
\left(\dfrac{m}{a}+1,\; \dfrac{b}{a}-1\right) \end{matrix};\; ax\right]
- \frac{\Gamma(m/a)}{\Gamma(m/a+1)}\right\}.
\label{eq:master_fox_wright_identity}
\end{equation}

For the case $m=1$ (ultra-radical), the identity takes a compact form:
\begin{equation}
M(1; a; b; x) = \frac{1}{a}\;
{}_1\Psi_1\!\left[\begin{matrix} \left(\dfrac{1}{a},\; \dfrac{b}{a}\right) \\[6pt]
\left(\dfrac{1}{a}+1,\; \dfrac{b}{a}-1\right) \end{matrix};\; ax\right].
\label{eq:ultraradical_fox_wright}
\end{equation}

\textbf{Verification of the first three terms.} Expanding the right-hand side of \eqref{eq:ultraradical_fox_wright} in powers of $x$ using $\Gamma(z+1) = z\,\Gamma(z)$:
\begin{align*}
n=0:&\quad \frac{1}{a}\cdot\frac{\Gamma(1/a)}{\Gamma(1/a+1)} = \frac{1}{a}\cdot\frac{1}{1/a} = 1,\\
n=1:&\quad \frac{1}{a}\cdot\frac{\Gamma((1+b)/a)}{\Gamma((1+b)/a)}\cdot\frac{ax}{1!} = x,\\
n=2:&\quad \frac{1}{a}\cdot\frac{\Gamma((1+2b)/a)}{\Gamma((1+2b)/a-1)}\cdot\frac{a^2x^2}{2!}
= \frac{1-a+2b}{2}\,x^2.
\end{align*}
These expressions coincide exactly with the first three terms of the master series
$M(1; a; b; x) = 1 + x + \frac{1-a+2b}{2}x^2 + \cdots$.

\textbf{Connection with the result of D\v{z}.~Belki\'c.} Formula (11.15) in \cite{Belkic2018} for an arbitrary power of the root of the trinomial equation $x = 1 + yx^{\alpha}$, expressed via ${}_1\Psi_1$, corresponds in the language of master series to identity \eqref{eq:ultraradical_fox_wright} at the fixed value $a=1$. Our formalism generalizes this result to arbitrary complex $a, b \in \mathbb{C}$ and, crucially, provides unambiguous branch identification both inside and beyond the radius of convergence via the geometric criterion.

\subsection{Connection with the Lambert--Tsallis Function}
\label{subsec:tsallis_connection}

The canonical identity $y = (1 + axy^b)^{1/a} = M(1;a;b;x)$ has a deep connection with the modern apparatus of deformed exponentials, actively applied in problems of plasma physics, control theory, and statistical mechanics. In particular, in the works of Brazilian authors \cite{SilvaRamos2019, Campos2026, Sousa2024}, the Lambert--Tsallis function (LTF) was proposed for the analytical solution of generalized trinomial equations $a_n x^n + a_m x^m + a_0 = 0$, based on the deformed exponential.

In the original formulation by Constantino Tsallis, the $q$-exponential is defined as
\begin{equation}
e_q(x) = \left[1 + (1-q)x\right]^{\frac{1}{1-q}}, \quad q \in \mathbb{R}.
\label{eq:q_exponential}
\end{equation}
As $q \to 1$, this function continuously transitions to the classical exponential: $\lim_{q \to 1} e_q(x) = e^x$.

In the works \cite{Campos2026, Sousa2024}, an alternative parametrization through $r = 1/(1-q)$ is used, in which the deformed exponential takes a more compact form
\begin{equation}
e_r(z) = \left(1 + \frac{z}{r}\right)^r, \quad r \in \mathbb{R}.
\label{eq:r_exponential}
\end{equation}
As directly noted in \cite{Campos2026}, the transition to the parameter $r$ ``simplifies algebraic manipulations and reveals important mathematical properties'' of the function, since in trinomial equations $r$ is directly equal to the ratio of the exponents: $r = (n-m)/m$.

Let us show that the canonical master equation of the ultraradical exactly coincides with the equation containing the Tsallis deformed exponential. Substituting the value $q = 1-a$ into formula \eqref{eq:q_exponential} and replacing the argument $x \to xy^b$, we obtain
\begin{equation}
e_{1-a}(xy^b) = \left[1 + (1 - (1-a))xy^b\right]^{\frac{1}{1-(1-a)}} = (1 + axy^b)^{\frac{1}{a}}.
\label{eq:canonical_tsallis_q}
\end{equation}
Thus, the canonical master equation of the ultraradical can be written in the form
\begin{equation}
y = e_{1-a}(xy^b) = M(1;a;b;x).
\label{eq:canonical_tsallis}
\end{equation}
In the parametrization of the Brazilian authors \cite{Campos2026, Sousa2024}, where $r = 1/(1-q) = 1/a$, the same equation takes the form
\begin{equation}
y = e_{1/a}(xy^b) = M(1;a;b;x).
\label{eq:canonical_tsallis_r}
\end{equation}

This identification has an important methodological consequence. Consider the limiting transition $a \to 0$. In the Tsallis parametrization, this corresponds to $q = 1-a \to 1$, and in the parametrization of \cite{Campos2026} --- to $r = 1/a \to \infty$. In both cases, the deformed exponential transitions to the ordinary one:
\begin{equation}
\lim_{a \to 0} e_{1-a}(xy^b) = \lim_{r \to \infty} e_r(xy^b) = e^{xy^b}.
\label{eq:tsallis_limit}
\end{equation}
Consequently, the canonical master equation as $a \to 0$ smoothly transitions to the transcendental equation $y = e^{xy^b}$, the solution of which is expressed through the classical Lambert $W$ function (for $b = -1$) or its generalizations. This demonstrates that the ultraradical $M(1;a;b;x)$ does not merely ``coexist'' with the Lambert--Tsallis function, but is its natural algebraic foundation: for $a \neq 0$ it describes the deformed exponential regime, and for $a \to 0$ --- the classical exponential regime.

In contrast to the approach of \cite{Campos2026}, where the Lambert--Tsallis function $W_r(z)$ is introduced as a ``black box'' for obtaining closed-form formulas through numerical solution of auxiliary polynomials, the master series theory provides an explicit parametric power series expansion with unambiguous identification of branches through the parameter $n$ and a strict geometric criterion for analytic continuation beyond the radius of convergence.

\section{Multivaluedness of solutions and branches of the complex logarithm}
\label{sec:branches}

\subsection{Multivaluedness of the canonical equation}

The master series $M(1; a; b; x)$, defined in Section~\ref{sec:master_series}, is a single-valued analytic function. However, the canonical equation
\begin{equation}
y^a = 1 + axy^b,
\label{eq:canonical_branches}
\end{equation}
considered in the complex domain for arbitrary $a, b \in \mathbb{C}$, has \emph{infinitely many} solutions, which is due to the multivaluedness of the operation of raising to a complex power.

Indeed, from equation \eqref{eq:canonical_branches} it formally follows:
\begin{equation}
y = (1 + axy^b)^{1/a}.
\label{eq:multivalued}
\end{equation}
The operation $z^{1/a} = \exp\left(\frac{1}{a}\ln z\right)$ is multivalued due to the multivaluedness of the complex logarithm:
\begin{equation}
\ln z = \ln|z| + i(\arg z + 2\pi k), \quad k \in \mathbb{Z}.
\label{eq:complex_log}
\end{equation}

\subsection{Method of branches: parametrization of solutions}

To systematically describe all solutions of equation \eqref{eq:canonical_branches}, we introduce a parametrization via the $a$-th roots of unity. Let $v$ be an arbitrary root of the equation $v^a = 1$, that is:
\begin{equation}
v_n = e^{f_n}, \quad f_n = \frac{2\pi i n}{a}, \quad n \in \mathbb{Z}.
\label{eq:vn_definition}
\end{equation}

\begin{pred}
\label{prop:branch_formula}
The family of solutions of equation \eqref{eq:canonical_branches} can be parametrized as follows (provided the master series converges):
\begin{equation}
y_n = v_n \cdot M(1; a; b; x V_n), \quad V_n = e^{b f_n}, \quad n \in \mathbb{Z},
\label{eq:branch_solution}
\end{equation}
where $M(1; a; b; \cdot)$ is the master series defined in Section~\ref{sec:master_series}.
\end{pred}

\textit{Proof.} Substitute $y_n$ from \eqref{eq:branch_solution} into equation \eqref{eq:canonical_branches}. Denote $Z_n = x V_n$ as the argument of the master series. By the definition of the master series, $Y_n = M(1; a; b; Z_n)$ satisfies the equation:
\begin{equation}
Y_n^a = 1 + a Z_n Y_n^b.
\label{eq:master_equation}
\end{equation}
Then:
\begin{align}
y_n^a &= (v_n Y_n)^a = v_n^a \cdot Y_n^a = 1 \cdot Y_n^a = 1 + a Z_n Y_n^b, \\
1 + axy_n^b &= 1 + ax(v_n Y_n)^b = 1 + ax \cdot v_n^b \cdot Y_n^b = 1 + a(x V_n) Y_n^b = 1 + a Z_n Y_n^b.
\end{align}
Comparing, we obtain $y_n^a = 1 + axy_n^b$, which completes the proof. $\square$

\begin{remark}
Formula \eqref{eq:branch_solution} shows that each solution $y_n$ is obtained by multiplying the single-valued master series by the factor $v_n = e^{2\pi i n / a}$, which encodes the branch of the complex logarithm. Thus, the multivaluedness of the solution is completely transferred to the parameter $n \in \mathbb{Z}$.
\end{remark}

\section{Other Canonical Master Equations}
\label{app:canonical}

\subsection{Special Cases with Zero Parameters}

If in the identity $y = \left( 1 + a x y^{b} \right)^{\frac{1}{a}} = M(1;a;b;x)$ we set $b=0$, we obtain the power series for the root of degree $a$:
\[
y = (1 + a x)^{\frac{1}{a}} = M(1;a;0;x) = 1 + x + (1 - a)\frac{x^{2}}{2} + (1 - a)(1 - 2a)\frac{x^{3}}{3!} + \cdots
\]

If $a=0$, we obtain the power series for the root of another equation:
\[
\lim_{a \rightarrow 0} y = {\lim_{a \rightarrow 0}\left( 1 + a x y^{b} \right)}^{\frac{1}{a}} = \lim_{a \rightarrow 0} M(1;a;b;x)
\]
\[
y = e^{x y^{b}} = M(1;0;b;x) = 1 + x + (1 + 2b)\frac{x^{2}}{2} + (1 + 3b)(1 + 3b)\frac{x^{3}}{3!} + \cdots
\]

For $a=0$, $b=0$:
\[
y = e^{x} = M(1;0;0;x) = 1 + x + \frac{x^{2}}{2} + \frac{x^{3}}{3!} + \cdots
\]

\subsection{Logarithmic Forms}

Replace the unknown $y$ with $e^{y}$ and take logarithm of both sides, setting the first parameter $m=0$. We obtain power series for roots of four more equations:

\begin{align*}
y &= x = M(0;0;0;x) = x \\
y &= x e^{b y} = M(0;0;b;x) = x + (2b)\frac{x^{2}}{2} + (3b)(3b)\frac{x^{3}}{3!} + \cdots \\
y &= \frac{\ln(1 + a x)}{a} = M(0;a;0;x) = x + ( - a)\frac{x^{2}}{2} + ( - a)( - 2a)\frac{x^{3}}{3!} + \cdots \\
y &= \frac{\ln\left( 1 + a x e^{b y} \right)}{a} = M(0;a;b;x) = x + (2b - a)\frac{x^{2}}{2} + (3b - a)(3b - 2a)\frac{x^{3}}{3!} + \cdots
\end{align*}
\subsection{Classification Table}
We organize three dimensions (three parameters) into a two-dimensional table, using column parity as the third dimension.

\begin{table}[ht]
\centering
\caption{Short Table of Canonical Master Equations}
\begin{tabular}{|c|c|c|c|c|}
\hline
$m = 0$ & $y = x$ & $y = x e^{b y}$ & $y = \frac{\ln(1 + a x)}{a}$ & $y = \frac{\ln\left( 1 + a x e^{b y} \right)}{a}$ \\
\hline
$m = 1$ & $y = e^{x}$ & $y = e^{x y^{b}}$ & $y = (1 + a x)^{\frac{1}{a}}$ & $y = \left( 1 + a x y^{b} \right)^{\frac{1}{a}}$ \\
\hline
 &  & $b \neq 0$ & $a \neq 0$ & $a \neq 0, b \neq 0$ \\
\hline
\end{tabular}
\end{table}
All these canonical equations are solved by a single unified master series.
\section{Master Series Identities}
\label{app:identities}

\subsection{Fundamental Identities}

\begin{enumerate}
\item \textbf{Exponential identity}:
\[
M(1;a;b;x) = \exp(M(0;a;b;x))
\]

\item \textbf{Scaling identity}:
\[
c \cdot M(0;a;b;x) = M(0;a/c;b/c;cx)
\]
Verified through power series:
\[
c \cdot M(0;a;b;x) = cx + (2b/c - a/c)\frac{(cx)^{2}}{2} + (3b/c - a/c)(3b/c - 2a/c)\frac{(cx)^{3}}{3!} + \cdots
\]

\item \textbf{Power identity}:
From the second identity through the first, we obtain the power identity:
\[
M^{c}(1;a;b;x) = M(1;a/c;b/c;cx)
\]

\item \textbf{Differential identity}:
Obtained by comparing derivatives of logarithm and power function:
\[
\frac{dM(1;a;b;x)}{dx} = M(1;a;b;x)\frac{dM(0;a;b;x)}{dx}
\]
\[
M(1;a;b;x) = \frac{dM(1;a;b;x)}{d\ln{M(1;a;b;x)}} = \frac{de^{M(0;a;b;x)}}{dM(0;a;b;x)} = e^{M(0;a;b;x)}
\]

\item \textbf{Isoroot identity}:
Master series inherit from master numbers:
\[
M(m;a;b;x) = M(m; - a;b - a;x)
\]

\item \textbf{Merge identity}:
From various transformations of master equations, for example:
\[
y = ve^{xy^{b}} = v \cdot M\left( 1;0;b;xv^{b} \right)
\]
\[
y = e^{xy^{b} + \ln v} = M\left( 1;0;b,0;x,\ln v \right)
\]
Therefore:
\[
M\left( 1;0;b,0;x,\ln v \right) = v \cdot M\left( 1;0;b;xv^{b} \right)
\]

\item \textbf{Arcsine identity}:
\[
M(0;2;1;x) = x + (3 - 2)(3 - 4)\frac{x^{3}}{3!} + (5 - 2)(5 - 4)(5 - 6)(5 - 8)\frac{x^{5}}{5!} + \cdots
\]
\[
M(0;2;1; - x) = - x - (3 - 2)(3 - 4)\frac{x^{3}}{3!} - (5 - 2)(5 - 4)(5 - 6)(5 - 8)\frac{x^{5}}{5!} + \cdots
\]
\[
M(0;2;1;x) = - M(0;2;1; - x)
\]

\item \textbf{Quadratic roots identity}:
\[
M(1;2;1;x) = \frac{1}{M(1;2;1; - x)}
\]
\end{enumerate}

\section{General Form of Master Equations}
\label{app:general_forms}

Table 2 shows canonical types of master equations. Any general form can be reduced to canonical form.

\subsection{Three-term exponential-type equation}

For example, three-term equation $pu = qe^{zu^b}$:
\[
pu = qe^{zu^{b}}
\]
Make substitution $u = y\frac{q}{p}$, $z = x\left( \frac{p}{q} \right)^{b}$, obtain canonical form:
\[
y = e^{xy^{b}}
\]

\subsection{Three-term algebraic equation}

Almost the same substitution reduces three-term equation $pu^{a} = q + zu^{b}$ to canonical form $y^{a} = 1 + xy^{b}$:
\[
pu^{a} = q + zu^{b}
\]
Make substitution:
\[
u = y\left( \frac{q}{p} \right)^{\frac{1}{a}},\quad z = qx\left( \frac{p}{q} \right)^{\frac{b}{a}}
\]
\[
qy^{a} = q + qxy^{b}
\]
\[
y^{a} = 1 + xy^{b}
\]

Completely different types of polynomials have practically identical solution through master series. The only difference is that $\left( \frac{q}{p} \right)^{\frac{1}{a}}$ has $a$ roots.

Root of the exponential-type equation $pu = qe^{zu^{b}}$:
\[
y = ve^{xy^{b}} = v \cdot M(1; 0; b; z),\quad z = xv^{b}
\]

Roots of the algebraic equation $py^{a} = q + xy^{b}$:
\[
py^{a} = q + xy^{b}
\]
\[
y = v \cdot M(1; a; b; z),\quad z = \frac{xv^{b}}{a q}
\]
\[
v = e^{\frac{\ln\left| \frac{q}{p} \right| + \left( \arg\left( \frac{q}{p} \right) + 2\pi n \right)i}{a}},\quad n \in \mathbb{Z}
\]

\subsection{Transformation to Logarithmic Forms}

Equations with $m=1$ can be converted to the $m=0$ form via the exponential substitution $u = e^{y}$ (introducing a new variable $y$):
\[
pu^{a} = q + xu^{b}
\]
\[
u = e^{y}
\]
\[
pe^{a y} = q + xe^{b y}
\]

The hyperbolic arcsine equation $e^{2y} = 1 + 2xe^{y}$
arises from the quadratic equation $y^{2} = 1 + 2xy$ via the substitution $y = e^{u}$.
\[
e^{2y} = 1 + 2xe^{y}
\]
Its roots are logarithms (meaning infinite branching) of the roots of the quadratic equation $y^{2} = 1 + 2xy$:
\[
y^{2} = 1 + 2xy
\]
\[
\operatorname{arsinh} x = \ln y
\]

The equation $u = e^{xu^{b}}$ can be reduced to logarithmic form by the substitution $u = e^y$:
\[
u = e^{xu^{b}}
\]
\[
u = e^{y}
\]
\[
e^{y} = e^{xe^{b y}}
\]
\[
y = xe^{b y}
\]

For $b=-1$, obtain Lambert W function equation:
\[
y = xe^{- y} = M(0;0; -1;x)
\]

\section{Convergence Radius of Master Series}
\label{app:convergence}
The master series $M(m; a; b; x) = m + \sum_{\ell=1}^{\infty} k_\ell x^\ell$, whose coefficients are given by
\begin{equation}
k_\ell = \frac{N_\ell}{\ell!} = \frac{1}{\ell!} \prod_{\gamma=1}^{\ell-1} (m - a\gamma + b\ell),
\label{eq:k_ell_def}
\end{equation}
has a finite radius of convergence $R$. According to the Cauchy--Hadamard formula:
\begin{equation}
\frac{1}{R} = \limsup_{\ell \to \infty} |k_\ell|^{\frac{1}{\ell}}.
\label{eq:cauchy_hadamard}
\end{equation}

To find the asymptotic behaviour as $\ell \to \infty$, we take the logarithm of $|k_\ell|$. The dominant contribution comes from the product $\prod_{\gamma=1}^{\ell-1}(b\ell - a\gamma)$, since $m$ is an additive constant. Factoring out $\ell$ from the product and replacing the sum of logarithms by an integral, using Stirling's formula $\ln(\ell!) \sim \ell\ln \ell - \ell$, we obtain:
\begin{equation}
\ln|k_\ell| \sim \ell \left( \int_{0}^{1} \ln|b - a \tau|\,d\tau + 1 \right) + o(\ell).
\label{eq:asymptotic_log}
\end{equation}

The integral can be evaluated analytically via the principal branch of the complex logarithm:
\begin{equation}
\int_{0}^{1} \ln|b - a \tau|\,d\tau = \operatorname{Re}\left( \frac{b\ln b - (b - a)\ln(b - a)}{a} \right) - 1.
\label{eq:integral_J}
\end{equation}

Consequently, the radius of convergence is:
\begin{equation}
R = \left| b^{-\frac{b}{a}} \cdot (b - a)^{\frac{b - a}{a}} \right|, \qquad a \neq 0,\; b \neq 0,\; b \neq a.
\label{eq:radius_final}
\end{equation}

For real parameters $a, b \in \mathbb{R}$ this formula simplifies to a compact form:
\begin{equation}
R = \frac{|1 - a/b|^{b/a}}{|b - a|}.
\label{eq:radius_compact}
\end{equation}

For degenerate cases the formula simplifies, and the case $b=a$ is elegantly resolved via the master identities:
\begin{itemize}
    \item \textbf{When $b = 0$:} $R = \frac{1}{|a|}$ (the classical binomial series).
    \item \textbf{When $a \to 0$:} $R = \frac{1}{|b| e}$ (the generalised exponential / Lambert function).
    \item \textbf{When $b = a$:} The equation takes the form $y^a = 1 + axy^a$, which is algebraically equivalent to $y = (1 - ax)^{-1/a}$. By the master identity $M(m;a;b;x) = M(m;-a;b-a;x)$, this corresponds precisely to the case $b=0$ with the parameter substitution $a \to -a$. Therefore, the radius of convergence is $R = \frac{1}{|a|}$.
\end{itemize}
\subsection{Behavior at the Convergence Boundary}
The behavior of power series at the boundary $|x| = R$ depends on the specific equation. For example:

\begin{itemize}
\item The series for $\sqrt{1+x}$ converges for $|x| \leq 1$
\item The series for $\frac{1}{\sqrt{1+x}}$ converges for $-1 < x \leq 1$
\end{itemize}

In this work, for brevity we use a single convergence condition $|x| < R$. It should be understood that the actual behavior depends on the specific equation:

\begin{itemize}
\item For some equations: $x < R$ 
\item For other equations: $x \leq R$
\end{itemize}

Derivatives and integrals of master series possess the same radius of convergence $R$ as the original series; the only possible difference lies in whether the series converges strictly inside the disk $|x| < R$ or also admits convergence on the boundary $|x| = R$.
\section{Isoroot Master Equations}
\label{app:isoroot}
Two different master equations are called \textbf{isoroot} if they have identical solution sets.

For example, if master equation $y^{a} = 1 + a x$, $y = M(1;a;0;x)$ is divided by $y^{a}$, we obtain master equation  $y^{- a} = 1 - a x y^{- a}$, $y = M(1; - a; - a;x)$. According to identity $M(1;a;b;x) = M(1; - a;b - a;x)$, we obtain identical roots in both cases.
\section{Analytical Continuations (Conjugate Master Series)}
\label{app:analytic_continuation}

Convergence radius of exponentials, sines and cosines is infinite. Binomial series has limited number of terms, since all subsequent master series terms equal zero. Therefore convergence question for this series is irrelevant.
\[
(1 + x)^{a} = M(1;1/a;0;ax)\quad ,\quad a \in \mathbb{N}
\]
\[
M(1;1/2;0;2x) = 1 + 2x + \frac{4x^{2}}{2}\left( 1 - \frac{1}{2} \right) + \frac{8x^{3}}{3!}\left( 1 - \frac{1}{2} \right)\left( 1 - \frac{2}{2} \right) = 1 + 2x + x^{2}
\]
\[
M(1;1/3;0;3x) = 1 + 3x + \frac{9x^{2}}{2}\left( 1 - \frac{1}{3} \right) + \frac{27x^{3}}{3!}\left( 1 - \frac{1}{3} \right)\left( 1 - \frac{2}{3} \right) + 0 + 0 + \cdots
\]

For other functions, question of convergence radius limitation is important. When one power series converges in region where original series diverges, and coincides with original function in common convergence region - this is called analytical continuation.

Master series $M(m;a;b;x)$ of some function converges for $|x| < R$, where $R$ is given by formula~\eqref{eq_radius}. For $|x| > R$ the original series diverges, and the function is represented by analytic continuation as an alternative master series $M(m';a';b';x')$ with different parameters and independent argument, converging in the required region.

Two or more master series representing the same equation solution in different regions of complex plane are called \textbf{conjugate}.

Consider equation of root of any degree $a$, including complex, of any number $x$:
\[
y = x^{\frac{1}{a}}
\]
$y$ determined through modulus $x$ and has multiple branches:
\[
x^{\frac{1}{a}} = |x|^{\frac{1}{a}}e^{\frac{\left( \arg(x) + 2\pi n \right)i}{a}},\ n \in \mathbb{Z}
\]

From modulus $x$, that is in positive number region, we can obtain root of any (including complex) degree $a$ through master series $|x|^{\frac{1}{a}} = M\left( 1;a;0;\frac{|x| - 1}{a} \right)$, but only if $|x| \leq 2$. If $|x| > 1$ need take $x$ itself to power $-1$, and root degree $a$ multiply by $-1$, to obtain conjugate master series converging in given complex region ($|x| \geq 1$).
\[
\left| x^{- 1} \right|^{\frac{1}{- a}} = M\left( 1; - a;0;\frac{\left| x^{- 1} \right| - 1}{- a} \right)
\]
\[
y = x^{\frac{1}{a}} = M\left( 1;ha;0;\frac{\left| x^{h} \right| - 1}{ha} \right)e^{\frac{\left( \arg(x) + 2\pi n \right)i}{a}},\ n \in \mathbb{Z},\ \begin{cases}
h = 1, & |x| < 1 \\
h = - 1, & |x| \geq 1
\end{cases}
\]

The hyperbolic arcsine equation has 2 roots, each with infinite number of branches. Each root determined by one of conjugate master series, depending on independent argument value. Transform this m=0 equation to m=1 equation:
\[
e^{2u} = 1 + 2we^{u}
\]
\[
u = \ln y
\]
\[
y^{2} = 1 + 2wy
\]
\[
y^{2} - 2wy - 1 = 0
\]
Use method for solving ABC\label{method_ABC}:
\[
Ay^{a} + By^{b} + C = 0
\]
\[
\operatorname{Re}(a) > \operatorname{Re}(b) > 0
\]
\textbf{Crucial condition:} the inequality $\operatorname{Re}(a) > \operatorname{Re}(b) > 0$ guarantees the correctness of the AB, BC, CA transformations and the proper selection of the conjugate series.
Any algebrao-transcendental trinomial $A Y^{a} + B Y^{b} + C = 0$ can be reduced to the condition $\operatorname{Re}(a) > \operatorname{Re}(b) > 0$ by multiplication with an appropriate power $Y^{f}$.

There exist 6 ways to transform this equation into $pY^{\alpha}=q+XY^{\beta}$. 2 ways use only permutation of equation terms. 4 ways, highlighted in yellow, use division of the entire equation by the unknown raised to a certain power.

\begin{table}[ht]
\centering
\caption{Isoroot Pairs of Master Equations for three-term algebraic equations}
\begin{tabular}{|c|c|}
\hline
$Ay^{a} = - C - By^{b}\quad (1a \equiv AB)$ & \cellcolor{yellow!20} $Cy^{- a} = - A - By^{b - a}\quad (1b)$ \\
\hline
\cellcolor{yellow!20} $By^{b - a} = - A - Cy^{- a}\quad (2a \equiv BC)$ & \cellcolor{yellow!20} $Ay^{a - b} = - B - Cy^{- b}\quad (2b)$ \\
\hline
\cellcolor{yellow!20} $Cy^{- b} = - B - Ay^{a - b}\quad (3a \equiv CA)$ & $By^{b} = - C - Ay^{a}\quad (3b)$ \\
\hline
\end{tabular}
\end{table}
\textbf{Remark.} The notation CA indicates that terms $C$ and $A$ are kept with $Y$, while the equation is divided by $Y$ raised to the power corresponding to the omitted term in the notation. For instance, the BC transformation results from dividing $AY^{a} + BY^{b} + C = 0$ by $Y^{a}$.

For solving three-term equation, three master equations suffice - remaining three are their mirror reflections and give same roots.
These transformations are required for more rigorous proofs and for deeper analysis of the ultra-radical. In programming, it is more convenient to use substitution lists\label{substitution_lists}:

\begin{itemize}
    \item $(1a \equiv Y(n))$: $p = A$, $X = -B$, $q = -C$, $\alpha = a$, $\beta = b$
    \item $(2a \equiv Y(h))$: $p = B$, $X = -C$, $q = -A$, $\alpha = b - a$, $\beta = -a$
    \item $(3a \equiv Y(k))$: $p = C$, $X = -A$, $q = -B$, $\alpha = -b$, $\beta = a - b$
\end{itemize}

All six transformations share a common solution method, called the pq-method\label{method_pq}:
\[
pY^{\alpha} = q + XY^{\beta}
\]
Define the parameters:
\[
f = \frac{\ln|q/p| + \bigl(\arg(q/p) + 2\pi N\bigr)i}{\alpha}, \quad
v = e^{f}, \quad
V = e^{\beta f}, \quad
Z = \frac{XV}{\alpha q}
\]
where $N \in \mathbb{Z}$. The solution is then expressed via the master series:
\[
Y(N) = v \cdot M(1;\alpha;\beta;Z).
\]
The integer parameter $N$ adopts different notations depending on the transformation:
\begin{itemize}
    \item For transformation AB: $N = n$ ($n \in \mathbb{Z}$)
    \item For transformation BC: $N = h$ ($h \in \mathbb{Z}$)
    \item For transformation CA: $N = k$ ($k \in \mathbb{Z}$)
\end{itemize}

All integer values $N$ for which the corresponding series converges constitute the complete solution set of the original equation. For real rational exponents, the roots generated by different $N$ repeat periodically; consequently, the equation possesses only finitely many solutions. If the exponents are complex, convergence for each particular $N$ and each transformation must be verified by the condition $|Z| < R$, where the convergence radius $R$ is given by formula (\ref{eq_radius}) in Section~\ref{app:convergence}.

For irrational or complex exponents, we recommend restricting attention to the principal root corresponding to $N=0$. For real exponents, convergence of the series does not depend on the value of $N$ and is determined solely by the transformation. If, moreover, $a > b > 0$, the choice of transformation is governed by the quantity
\[
T = \left| \frac{b}{A} \right|^{b} {\left| \frac{B}{a} \right|^{a} \left| \frac{a - b}{C} \right|}^{a - b}.
\]

When $T < 1$, all roots of the equation are described by transformation (1a) (or the equivalent (1b)). If $T \ge 1$, then $a^{\prime}-b^{\prime}$ roots are given by transformation (2a) and the remaining $b^{\prime}$ roots by transformation (3a). Here $a^{\prime}$ and $b^{\prime}$ denote the numerators of the rational fractions $a$ and $b$ reduced to a common denominator (i.e., the integers obtained after reduction).

Using the ABC transformations (or, equivalently, the $n$, $h$, $k$ substitutions), one can compute the values of the ultra‑radical $\sqrt[{n;a;b}]{x}$, i.e. the solutions of $y^{a}=1+a x y^{b}$, even outside the convergence radius of the original power series. The crucial point is that the ultra‑radical performs a strict conjugation among different power series: when the original series diverges ($|x| \ge R$), the corresponding branch is continued analytically by one of the two conjugate series (the $h$-series or the $k$-series). The choice of the correct conjugate series for each branch $n$ is determined unambiguously by the geometric criterion described in Section~\ref{sec_geometric_criterion} (see the conjugation rule).

Although the ABC method is in principle applicable to the equation $AY^{r}+BY^{s}+C=0$, the geometric sector‑selection criterion described in Section~\ref{sec_geometric_criterion} was originally devised for the canonical form of the ultra‑radical. For fractional or irrational exponents, a direct application of this criterion within the ABC transformations would require an extension of the algorithm. A simpler and more reliable route is to reduce the original equation to canonical form, solve it via the ultra‑radical, and then return to the original variables.

The required transformations are:
\begin{align}
y^{r} &= 1 + r x y^{s} \quad \text{(canonical form of the ultra‑radical)} \label{eq:canonical_en} \\
x &= \frac{B}{r C}\left(-\frac{C}{A}\right)^{\!s/r} \quad \text{(principal branch)} \label{eq:x_transform} \\
Y_n &= y_n \left(-\frac{C}{A}\right)^{\!1/r} 
     = \left(-\frac{C}{A}\right)_{0}^{\!1/r} \cdot U_{n}(r;s;x_0) \label{eq:Y_transform}
\end{align}
where $U_{n}(r;s;x_0) = \sqrt[{n;r;s}]{x_0}$ is the ultra‑radical evaluated at the argument $x_0$ obtained from formula~\eqref{eq:x_transform} using the principal branch ($n=0$) of the complex power.

Hence, the complete solution of the original three‑term equation is obtained in three steps:
\begin{enumerate}
    \item Compute $x_0$ via formula~\eqref{eq:x_transform} (principal branch).
    \item Find the canonical roots $y_n = U_{n}(r;s;x_0)$ through the ultra‑radical.
    \item Recover the desired roots $Y_n$ using formula~\eqref{eq:Y_transform} (again the principal branch).
\end{enumerate}
This approach ensures that branch conjugation always respects the geometric criterion inherent in the definition of the ultra‑radical.

\paragraph{Conjugate Series for the Hyperbolic Arcsine}
We now return to the problem of obtaining all power‑series representations and their analytic continuations for the hyperbolic arcsine, derived directly from its defining equation.

For $T < 1$, both quadratic equation roots determined through master equation:
\[
Ay^{2} = - C - By
\]
Solved by master method:
\[
py^{r} = q + xy^{s}
\]
\[
y = v \cdot M(1;r;s;z),\quad z = \frac{xv^{s}}{rq},\quad v = e^{\frac{\ln\left| \frac{q}{p} \right| + \left( \arg\left( \frac{q}{p} \right) + 2\pi n \right)i}{r}},\ n \in \mathbb{Z}
\]

$p=A$, $r=2$, $q=-C$, $x=-B$, $s=1$:
\[
y = v M(1;2;1;z),\quad z = Bv/2C,\quad v = \frac{\ln\left| \frac{C}{A} \right| + \left( \arg\left( \frac{- C}{A} \right) + 2\pi n \right)i}{2},\ n \in \mathbb{Z}
\]

For $n=0$, $v = 1$, $z = w$

For $n=1$, $v = - 1$, $z = - w$

Use identity:
\[
M(1;2;1;w) = \frac{1}{M(1;2;1; - w)}
\]

For $n=0$:
\[
y_{0} = M(1;2;1;w)
\]

For $n=1$:
\[
y_{1} = - \frac{1}{M(1;2;1;w)} = - \frac{1}{y_{0}}
\]

For $|w| \geq 1$, power series of one root converges under transformation (2) $By^{b - a} = - A - Cy^{- a}$, power series of other root converges under transformation (3) $Cy^{- b} = - B - Ay^{a - b}$.

Under transformation (2):
\[
By^{b - a} = - A - Cy^{- a}
\]
\[
2wy^{- 1} = 1 - y^{- 2}
\]

Use method y(h):
\[
py^{r} = q + xy^{s}
\]
\[
y = v \cdot M(1;r;s;z),\quad z = \frac{xv^{s}}{rq},\quad v = e^{\frac{\ln\left| \frac{q}{p} \right| + \left( \arg\left( \frac{q}{p} \right) + 2\pi h \right)i}{r}},\ h \in \mathbb{Z}
\]
\[
y = v \cdot M(1; - 1; - 2;z),\quad z = v^{- 2},\quad v = 2w
\]

Under transformation (3):
\[
Cy^{- b} = - B - Ay^{a - b}
\]
\[
y^{- 1} = - 2w + y^{1}
\]

Use method y(k):
\[
py^{r} = q + xy^{s}
\]
\[
y = v \cdot M(1;r;s;z),\quad z = \frac{xv^{s}}{rq},\quad v = e^{\frac{\ln\left| \frac{q}{p} \right| + \left( \arg\left( \frac{q}{p} \right) + 2\pi k \right)i}{r}},\ k \in \mathbb{Z}
\]
\[
y = v \cdot M(1; - 1;1;z),\quad z = \frac{v}{2w},\quad v = - \frac{1}{2w}
\]

We obtained analytical continuation of two roots of equation $y^{2} - 2wy - 1 = 0$, when $|w| \geq 1$:
\[
y(h) = 2w \cdot M\left( 1; - 1; - 2;\frac{1}{4w^{2}} \right)
\]
and
\[
y(k) = - \frac{1}{2w} \cdot M\left( 1; - 1;1; - \frac{1}{4w^{2}} \right)
\]

Total: we obtained 2 pairs of power series for $|w| < 1$:
\[
y_{0} = M(1;2;1;w)
\]
\[
y_{1} = - M(1;2;1; - w)
\]
and 2 pairs of power series for $|w| \geq 1$:
\[
y(h) = 2w \cdot M\left( 1; - 1; - 2;\frac{1}{4w^{2}} \right)
\]
\[
y(k) = - \frac{1}{2w} \cdot M\left( 1; - 1;1; - \frac{1}{4w^{2}} \right)
\]
which are roots of quadratic equation:
\[
y^{2} - 2wy - 1 = 0
\]

Choice of roots $h$ and $k$ performed by sector method (see section 2).

Thus, complete equation solution is not one series, but network of conjugate series, connected by transformations and covering entire complex plane for each branch of multi-valued function.

Hyperbolic arcsine expansion:
\[
\ln M(1;r;s;w) = M(0;r;s;w)
\]
\[
\operatorname{arsinh}_{0}(w) = \ln\left( y_{0} \right) = M(0;2;1;w) + 2\pi ki,\ k \in \mathbb{Z}
\]
\[
M(0;2;1;w) = w + (3 - 2)(3 - 4)\frac{w^{3}}{3!} + (5 - 2)(5 - 4)(5 - 6)(5 - 8)\frac{w^{5}}{5!} + \cdots
\]
\[
\ln M\left( 1; - 1; - 2;\frac{1}{4w^{2}} \right) = M\left( 0; - 1; - 2;\frac{1}{4w^{2}} \right)
\]
\[
\operatorname{arsinh}_{h}(w) = \ln\left( y_{0} \right) = \ln(2w) + M\left( 0; - 1; - 2;\frac{1}{4w^{2}} \right)
\]

Other expansions:
\[
y_{1} = - M(1;2;1; - w)
\]
\[
\operatorname{arsinh}_{1}(w) = \ln\left( y_{1} \right) = M(0;2;1; - w) + (\pi + 2\pi k)i,\ k \in \mathbb{Z}
\]
\[
M(0;2;1; - w) = - w - (3 - 2)(3 - 4)\frac{w^{3}}{3!} - (5 - 2)(5 - 4)(5 - 6)(5 - 8)\frac{w^{5}}{5!} + \cdots
\]

\section{Master Numbers (Generalization of Factorials and Powers)}
\label{app:master_numbers}
\textbf{Clarification:} Master numbers $N(m; a; b; \ell)$ are \textbf{not} power series — they are the \textbf{coefficients} in the power series expansion of master functions. Each master number is a \textbf{finite product} that can be computed directly without series expansion.

A master number of order $\ell$ is the product of factors $(m - a\gamma + b\ell)$. The number of factors equals $\ell-1$. Therefore for $\ell<2$ the master number returns 1. We denote master numbers inline (N), and use a short notation analogous to that employed for various factorial generalizations, powers, and other special numbers.
\[
N(m; a; b; \ell) = \mstr{m}{a}{b}{\ell} = \prod_{\gamma > 0}^{\gamma < \ell}(m - a\gamma + b\ell)
\]
\[
\mstr{m}{a}{b}{\ell} = (m - a + b\ell)(m - 2a + b\ell)\cdots\left( m - a(\ell - 1) + b\ell \right)
\]

Important identity:
\[
\mstr{m}{a}{b}{\ell} = (m - a + b\ell)(m - 2a + b\ell)\cdots\left( m - a(\ell - 2) + b\ell \right)\left( m - a(\ell - 1) + b\ell \right)
\]
\[
\mstr{m}{-a}{b-a}{\ell} = \left( m - a(\ell - 1) + b\ell \right)\left( m - a(\ell - 2) + b\ell \right)\cdots(m - 2a + b\ell)(m - a + b\ell)
\]

Conclusion:
\[
N(m; a; b; \ell) = N(m; -a; b-a; \ell)
\]
\[
\mstr{m}{a}{b}{\ell} = \mstr{m}{-a}{b-a}{\ell}
\]

Example:
\[
\mstr{1}{-2}{0}{4} = (1 + 2)(1 + 4)(1 + 6) = 3 \cdot 5 \cdot 7
\]
\[
\mstr{1}{2}{2}{4} = (1 + 8 - 2)(1 + 8 - 4)(1 + 8 - 6) = 7 \cdot 5 \cdot 3
\]

\begin{table}[ht]
\centering
\caption{Examples of Master Numbers Usage}
\begin{tabular}{|l|l|l|l|}
\hline
\textbf{Function} & \textbf{Notation} & \textbf{Master Number} & \textbf{Example} \\
\hline
Gamma function & $\Gamma(\ell)$ & $\mstr{0}{1}{1}{\ell}$ & $\mstr{0}{1}{1}{3} = \Gamma(3) = 2 \times 1$ \\
\hline
Factorial & $\ell!$ & $\mstr{1}{1}{1}{\ell}$ & $\mstr{1}{1}{1}{3} = 3! = 3 \times 2$ \\
\hline
$a$-fold factorial & ${n!}_{(a)}$ & $\mstr{n-a\ell}{a}{a}{\ell+1}$ & $\ell = \left\lceil \frac{n}{a} \right\rceil$ \\
\hline
Falling factorial & $(n)_{k}$ & $\mstr{n+1}{1}{0}{k+1}$ & $(5)_{3} = \mstr{6}{1}{0}{4} = 5 \times 4 \times 3$ \\
\hline
Binomial coefficient & $\binom{n}{k}$ & $\frac{\mstr{n+1}{1}{0}{k+1}}{k!}$ & $\binom{4}{2} = \frac{\mstr{5}{1}{0}{3}}{2!} = 6$ \\
\hline
Rising factorial & $n^{(k)}$ & $\mstr{n-1}{-1}{0}{k+1}$ & $3^{(3)} = \mstr{2}{-1}{0}{4} = 3 \times 4 \times 5$ \\
\hline
\end{tabular}
\end{table}
Examples of triple factorial from numbers 5, 6, 7:
\begin{align*}
5!!! &= {5!}_{(3)} = \mstr{5-3\ell}{3}{3}{\ell+1},\ \ell = \left\lceil \frac{5}{3} \right\rceil = 2 \\
&= \mstr{-1}{3}{3}{3} = (-1 - 3 + 9)(-1 - 6 + 9) = 5 \times 2
\end{align*}
\begin{align*}
{6!}_{(3)} &= \mstr{6-3\ell}{3}{3}{\ell+1},\ \ell = \left\lceil \frac{6}{3} \right\rceil = 2 \\
&= \mstr{0}{3}{3}{3} = (9 - 3)(9 - 6) = 6 \times 3
\end{align*}
\begin{align*}
{7!}_{(3)} &= \mstr{7-3\ell}{3}{3}{\ell+1},\ \ell = \left\lceil \frac{7}{3} \right\rceil = 3 \\
&= \mstr{-2}{3}{3}{4} = (-2 - 3 + 12)(-2 - 6 + 12)(-2 - 9 + 12) = 7 \times 4 \times 1
\end{align*}
\subsection{Merge Operation for Master Numbers}
Master numbers can be merged into composite structures:
\[
\mstr{m}{a}{b_1,b_2,\cdots}{\ell_1,\ell_2,\cdots} = \prod_{\gamma = 1}^{\gamma < \ell_1 + \ell_2 + \cdots}\left( m - a\gamma + b_1\ell_1 + b_2\ell_2 + \cdots \right)
\]
\[
\mstr{m}{a}{b_1,b_2,b_3}{\ell_1,\ell_2,\ell_3} = \prod_{\gamma = 1}^{\ell_1 + \ell_2 + \ell_3 - 1}\left( m - a\gamma + b_1\ell_1 + b_2\ell_2 + b_3\ell_3 \right)
\]

This operation allows only paired swaps of indices (e.g., $(b_1 , \ell_1)$ with $(b_2 , \ell_2)$) and generalizes multi-indexed series expansions.
\section{Function Merging Operation}
\label{app:function_merging}

The merge operation (denoted by symbol $@$) allows combining multiple series into a single unified structure. This operation is essential for solving equations with multiple terms, where each term requires its own master series with distinct independent arguments.

The merge operation for master numbers is defined as:
\[
\mstr{m}{a}{b_1}{\ell_1} @ \mstr{m}{a}{b_2}{\ell_2} @ \cdots = \mstr{m}{a}{b_1, b_2, \cdots}{\ell_1, \ell_2, \cdots} = \prod_{\gamma = 1}^{\ell_1 + \ell_2 + \cdots - 1} \left( m - a\gamma + b_1\ell_1 + b_2\ell_2 + \cdots \right)
\]

\subsection{Merging Power Series}

For power series ($m = 1$), the merge operation is performed as follows:
\[
\msr{1}{a}{b_1}{x_1} @ \msr{1}{a}{b_2}{x_2} @ \cdots = \msr{1}{a}{b_1, b_2, \cdots}{x_1, x_2, \cdots}
\]

Merging of power series is analogous to ordinary multiplication of power series, where each term of one series is multiplied by each term of all other series, with the exception of the master-number factors.

\subsection{Features of Term Merging}

\begin{enumerate}
\item \textbf{Independent fractions are multiplied}:
\[
\frac{x_1^{\ell_1}}{\ell_1!} \cdot \frac{x_2^{\ell_2}}{\ell_2!} \cdot \frac{x_3^{\ell_3}}{\ell_3!}
\]

\item \textbf{Master numbers form a unified monolith}: Master numbers of all merged terms form a single structure that cannot be obtained by simple multiplication of the original master numbers:
\[
\mstr{m}{a}{b_1, b_2, b_3}{\ell_1, \ell_2, \ell_3} = \prod_{\gamma = 1}^{\ell_1 + \ell_2 + \ell_3 - 1} \left( m - a\gamma + b_1\ell_1 + b_2\ell_2 + b_3\ell_3 \right)
\]

\item \textbf{Zero term handling}: For power series ($m=0$), the zero term is $w_0 = 0$. During merging, the zero term is temporarily set to $w_0 = 1$. After completing the merge operation, 1 is subtracted from the resulting monolith.

\item \textbf{Special case $\ell = 1$}: A master number for $\ell = 1$ always equals 1, regardless of parameters, since the number of master factors is one less than $\ell$. However, the monolith of the first terms no longer equals 1, because the sum of the ordinal numbers of these terms is greater than 1:
\[
\mstr{m}{a}{b_1}{1} = 1,\quad \mstr{m}{a}{b_2}{1} = 1,\quad \mstr{m}{a}{b_3}{1} = 1
\]
\[
\mstr{m}{a}{b_1}{1} @ \mstr{m}{a}{b_2}{1} = \mstr{m}{a}{b_1, b_2}{1, 1} = (m - a + b_1 + b_2)
\]
\[
\mstr{m}{a}{b_1}{1} @ \mstr{m}{a}{b_2}{1} @ \mstr{m}{a}{b_3}{1} = \mstr{m}{a}{b_1, b_2, b_3}{1, 1, 1} = (m - a + b_1 + b_2 + b_3)(m - 2a + b_1 + b_2 + b_3)
\]
\end{enumerate}

For brevity we introduce a master-series notation analogous to the master-numbers notation:
\[
M(m;a;b;x) = \msr{m}{a}{b}{x}
\]

The binomial series can be extended to any number of terms, but here too the merge operation is used.
\[
\left( 1 + x_1 + x_2 + \cdots \right)^{a} = \msr{1}{1/a}{0, 0, \cdots}{a x_1, a x_2, \cdots}\quad ,\quad a \in \mathbb{N}
\]
\[
\left( 1 + \frac{x}{2} \right)^{2} = \msr{1}{1/2}{0}{x} = 1 + x + \frac{x^{2}}{4}
\]
\[
\left( 1 + \frac{x}{3} \right)^{3} = \msr{1}{1/3}{0}{x} = 1 + x + \frac{x^{2}}{3} + \frac{x^{3}}{27}
\]
\[
\left( 1 + \frac{x_1}{2} + \frac{x_2}{2} \right)^{2} = \msr{1}{1/2}{0,0}{x_1,x_2} = 1 + x_1 + \frac{x_1^2}{4} + x_2 + \frac{x_1x_2}{2} + \frac{x_2^2}{4}
\]

This demonstrates that multinomial coefficients are essentially special cases of merged master numbers. The most important application of function merging is obtaining analytical solutions to polynomials with an arbitrary number of terms.

The most detailed study of the merge operation can be found in the JavaScript code \cite{Berezin2024Lab}.

Solving polynomials with an arbitrary number of terms (algebraic equations of any degree, including complex exponents):

Three-term → one master series

Four-term → merge two master series (two‑core series)

Five-term → merge three master series (three‑core series)
\[
py^{a} = q + x_{1}y^{b_{1}} + x_{2}y^{b_{2}} + \cdots = q \cdot \msr{1}{1}{b_{1}/a, b_{2}/a, \cdots}{\frac{x_{1}(q/p)^{b_{1}/a}}{q}, \frac{x_{2}(q/p)^{b_{2}/a}}{q}, \cdots}
\]
\[
py^{a} = q + x_{1}y^{b_{1}} + x_{2}y^{b_{2}} + \cdots = p\left( v \cdot \msr{1}{a}{b_{1}, b_{2}, \cdots}{\frac{x_{1}v^{b_{1}}}{a q}, \frac{x_{2}v^{b_{2}}}{a q}, \cdots} \right)^{a}
\]
\[
y = v \cdot \msr{1}{a}{b_{1}, b_{2}, \cdots}{\frac{x_{1}v^{b_{1}}}{a q}, \frac{x_{2}v^{b_{2}}}{a q}, \cdots},\quad v = \sqrt[a]{\left| \frac{q}{p} \right|} e^{\frac{\left( \arg\left( \frac{q}{p} \right) + 2\pi n \right)i}{a}},\ n \in \mathbb{Z}
\]
\section{Roots of Other Equation Types via Master Series}
\label{app:other_equations}

The master series provides a unified framework for representing solutions to various types of equations beyond the canonical forms. This demonstrates the universality of the approach and enables systematic analysis of diverse mathematical functions.

\subsection{Elementary and Special Functions}

\begin{align*}
e &= \msr{1}{0}{0}{1} \\
\exp(x) &= \msr{1}{0}{0}{x} \\
\cos(x) &= \overline{\overline{\msr{1}{0}{0}{ix}}} \\
\sin(x) &= \overline{\msr{1}{0}{0}{ix}} \\
\sinh(x) &= \overline{\msr{1}{0}{0}{x}} \\
\cosh(x) &= \overline{\overline{\msr{1}{0}{0}{x}}} \\
\pi i &= 4\overline{\msr{0}{1}{0}{i}} \\
\operatorname{artanh}(x) &= \overline{\msr{0}{1}{0}{x}} \\
W_0(x) &= \msr{0}{0}{-1}{x} \\
\msr{0}{0}{b}{x} &= \frac{W_0(-bx)}{-b} \\
\ln(1+x) &= \msr{0}{1}{0}{x}, \quad |x| < 1
\end{align*}

where the overline notation denotes:
\[
\overline{\msr{m}{a}{b}{x}} = \frac{\msr{m}{a}{b}{x} - \msr{m}{a}{b}{-x}}{2}, \quad
\overline{\overline{\msr{m}{a}{b}{x}}} = \frac{\msr{m}{a}{b}{x} + \msr{m}{a}{b}{-x}}{2}
\]

\subsection{Composite Master Series}

More complex functions can be expressed through composite master series:
\begin{align*}
\sec x &= \msr{1}{1}{1}{1 - \cos x} \\
\csc x &= \msr{1}{1}{1}{1 - \sin x} \\
\tan x &= \sin x \cdot \msr{1}{1}{1}{1 - \cos x} \\
\cot x &= \cos x \cdot \msr{1}{1}{1}{1 - \sin x}
\end{align*}

Although trigonometric series have infinite convergence radii, they appear here as independent arguments of geometric series. To obtain power series in $x$, substitute sine and cosine with their power series, expand brackets, and group terms by $x$.

\subsection{Equations with Trigonometric Dependencies}

Merging series provides solutions to other equation types. For example:
\[
y = x \cdot \cos(by) = \frac{x \cdot e^{byi} + x \cdot e^{-byi}}{2} = \msr{0}{0}{bi, -bi}{x/2, x/2}
\]

Similarly, for equations with Bessel functions or other special functions, appropriate master series representations can be derived through suitable transformations and merge operations.

\subsection{Integral Forms and Derivatives}

For denoting power series using master series integrals, we can use a fifth parameter. Its value indicates how many derivatives to take from the given series to obtain the master series with specified parameters.

\begin{align*}
M(m;a;b;x) &= \msr{m}{a}{b}{x} = m + x + \sum_{\ell = 2}^{\infty}\left( \frac{x^{\ell}}{\ell!}\prod_{\gamma = 1}^{\ell - 1}(m - a\gamma + b\ell) \right) \\
M(m;a;b;x;1) &= \msr[1]{m}{a}{b}{x} = mx + \frac{x^{2}}{2} + \sum_{\ell = 2}^{\infty}\left( \frac{x^{\ell + 1}}{(\ell + 1)!}\prod_{\gamma = 1}^{\ell - 1}(m - a\gamma + b\ell) \right) \\
&= mx + \frac{x^{2}}{2} + (m - a + 2b)\frac{x^{3}}{3!} + (m - a + 3b)(m - 2a + 3b)\frac{x^{4}}{4!} + \cdots \\
M(m;a;b;x;2) &= \msr[2]{m}{a}{b}{x} = m\frac{x^{2}}{2} + \frac{x^{3}}{3!} + \sum_{\ell = 2}^{\infty}\left( \frac{x^{\ell + 2}}{(\ell + 2)!}\prod_{\gamma = 1}^{\ell - 1}(m - a\gamma + b\ell) \right) \\
&= m\frac{x^{2}}{2} + \frac{x^{3}}{3!} + (m - a + 2b)\frac{x^{4}}{4!} + (m - a + 3b)(m - 2a + 3b)\frac{x^{5}}{5!} + \cdots
\end{align*}

\subsection{Systematic Classification Implications}

The ability to represent diverse functions through master series enables:

\begin{itemize}
\item \textbf{Unified convergence analysis}: All represented functions share common convergence criteria
\item \textbf{Systematic analytical continuation}: Continuation methods developed for master series apply to all represented functions
\item \textbf{Computational efficiency}: Single algorithm handles multiple function types
\item \textbf{Error analysis}: Uniform error estimation across different function classes
\item \textbf{Symbolic manipulation}: Consistent transformation rules for various equation types
\end{itemize}

This systematic approach offers a unified perspective on function analysis, moving from a collection of specialized methods toward a more integrated computational framework.

\section{Relationship Between Master Series and Classical Hypergeometric Functions}
\label{sec:master_vs_hypergeometric}

The master series $M(m; a; b; x)$ and the hypergeometric function ${}_{p}F_{q}$ are two distinct ways of parameterizing power series. In this section we show that in certain degenerate cases the two formalisms yield identical expansions for classical functions, thus establishing precise correspondences between their parameters.

\subsection{Definitions}

The master series is defined as
\[
M(m; a; b; x) = m + x + \sum_{\ell=2}^{\infty} \frac{x^{\ell}}{\ell!}
\prod_{\gamma=1}^{\ell-1} (m - a\gamma + b\ell),
\]
where $m, a, b \in \mathbb{C}$.

The hypergeometric series:
\[
{}_{p}F_{q}\!\left( a_{1},\dots ,a_{p}; b_{1},\dots ,b_{q}; x \right)
= \sum_{k=0}^{\infty} \frac{(a_{1})_{k}\cdots (a_{p})_{k}}{(b_{1})_{k}\cdots (b_{q})_{k}} 
\frac{x^{k}}{k!},
\]
where $(u)_{k}=u(u+1)\dots (u+k-1)$ is the Pochhammer symbol.

\subsection{Examples of Exact Correspondences}

\paragraph{Binomial series}
\[
(1+x)^{\alpha} = {}_{1}F_{0}(\alpha;;-x)
= M\!\left(1; \frac{1}{\alpha}; 0; \alpha x\right).
\]

\paragraph{Inverse power function}
\[
(1-x)^{-\alpha} = {}_{1}F_{0}(\alpha;;x)
= M\!\left(1; -\frac{1}{\alpha}; 0; \alpha x\right).
\]

\paragraph{Natural logarithm}
\[
\ln(1+x) = x\;{}_{2}F_{1}(1,1;2;-x)
= M(0;1;0;x).
\]

\paragraph{Exponential function}
\[
e^{x} = {}_{0}F_{0}(;;x) = {}_{1}F_{1}(1;1;x)
= M(1;0;0;x).
\]

\subsection{Interpretation of the Coincidences}

The equalities above reflect the fact that when the parameters \(a\) or \(b\) take special values (zero or related to the exponent \(\alpha\)), the master series \(M(m;a;b;x)\) reduces to elementary functions --- the binomial \((1+x)^{\alpha}\), the logarithm \(\ln(1+x)\), or the exponential \(e^x\). These elementary functions, in turn, admit classical hypergeometric representations.

It is crucial to note that this correspondence occurs only in these \emph{degenerate} parameter configurations. For genuine nonlinear master equations (when both \(a\neq0\) and \(b\neq0\)), the solutions --- ultra-radicals and ultra-logarithms --- are not hypergeometric functions, and their analytic continuation, branching structure, and stability properties are governed by the geometric criterion and merging operation developed in this work, which have no counterpart in classical hypergeometric theory.

\subsection{Conclusion}

The coincidence of master series with hypergeometric series in degenerate cases confirms that the master formalism correctly generalizes classical analysis. However, the complete theory of analytic continuation, branching, and stability of solutions to nonlinear master equations presented in this work is self‑contained and cannot be reduced to known results about hypergeometric functions.

\subsection*{The Problem of Direct Correspondence}

The classical Gauss hypergeometric function
\[
{}_2F_1(a,b;c;x)=1+\frac{ab}{c}x+\frac{a(a+1)b(b+1)}{c(c+1)}\frac{x^2}{2!}+\cdots
\]
uses \emph{Pochhammer symbols} \((a)_n=a(a+1)\dots(a+n-1)\).

The master series \(M(m;a;b;x)\) uses \emph{master numbers}
\[
N(m;a;b;\ell)=\prod_{\gamma=1}^{\ell-1}(m-a\gamma+b\ell).
\]

To establish an exact correspondence between these two types of coefficients one must pass to a more general structure --- the \textbf{hyper‑master}.

\subsection{Definition of the Hyper‑Master}

Let the following sets of parameters be given:
\begin{itemize}
    \item In the numerator: \((m_j, a_j, b_j)\) for \(j=1,\dots,p\)
    \item In the denominator: \((M_k, A_k, B_k)\) for \(k=1,\dots,q\)
\end{itemize}

The hyper‑master kernel (without constant term) is defined as
\[
_{p}^{q}H\bigl(\{m_j;a_j;b_j\};\{M_k;A_k;B_k\};x\bigr)
= x+\sum_{\ell=2}^{\infty}\frac{x^\ell}{\ell!}
\prod_{\gamma=1}^{\ell-1}
\frac{\prod_{j=1}^p (m_j-a_j\gamma+b_j\ell)}
     {\prod_{k=1}^q (M_k-A_k\gamma+B_k\ell)}.
\]

The hyper‑master with unit constant term is denoted by \(_{p}^{q}H_1(\cdots;x)\).

\subsection*{Exact Representation of \({}_2F_1\) via the Hyper‑Master}

The function \({}_2F_1(a,b;c;x)\) is expressed through the hyper‑master as
\[
{}_2F_1(a,b;c;x)
= {}_{1}^{-;2}H_1\Bigl(
   1;\; -\tfrac{1}{a},0;\; -\tfrac{1}{b},0;\; -\tfrac{1}{c},0;\;
   \tfrac{ab}{c}x\Bigr),
\]
where the symbol ``\(-\)'' in the upper index indicates that the parameter \(m\) is the same in all master numbers and equals one.

In expanded form this yields the series
\[
{}_2F_1(a,b;c;x)=1+\frac{ab}{c}x+
\sum_{\ell=2}^{\infty}
\frac{\bigl(\frac{ab}{c}x\bigr)^\ell}{\ell!}
\prod_{\gamma=1}^{\ell-1}
\frac{(1+\frac{1}{a}\gamma)(1+\frac{1}{b}\gamma)}
     {(1+\frac{1}{c}\gamma)},
\]
which, after elementary transformations, coincides with the classical hypergeometric series.

\subsection*{Key Features of the Correspondence}

\begin{enumerate}
    \item \textbf{Inversion of parameters}. The hypergeometric parameters \(a,b,c\) enter the hyper‑master as \(-\tfrac{1}{a},-\tfrac{1}{b},-\tfrac{1}{c}\). This stems from the fact that in the master equation \(y^a=1+axy^b\) the exponent \(a\) stands \emph{in the power of the unknown}, whereas in the hypergeometric function the parameters are \emph{arguments of the Pochhammer symbol}.

    \item \textbf{Scaling of the argument}. The argument \(x\) of the hypergeometric function in the hyper‑master representation is multiplied by the factor \(\frac{ab}{c}\), reflecting the canonical form of the master equation where the independent variable is always multiplied by the power‑law parameter.

    \item \textbf{Generality}. The hyper‑master with arbitrary sets of parameters \(\{a_j,b_j\}\), \(\{A_k,B_k\}\) naturally generalizes not only \({}_2F_1\) but the whole hierarchy \({}_pF_q\). At the same time it retains the algebraic structure that admits the merging operation (@), opening a way to solving multi‑term equations inaccessible to the classical hypergeometric apparatus.
\end{enumerate}

\subsection*{Open Directions}

\begin{itemize}
    \item \textbf{Differential equations}. What differential equation does a general hyper‑master \(_{p}^{q}H\) satisfy?
    \item \textbf{Merging operation}. What is the combinatorial meaning of merging hyper‑masters and how is it related to representing solutions of equations with several nonlinear terms?
    \item \textbf{Analytic continuation}. Can the geometric branch‑selection criterion developed for ultra‑radicals be extended to hyper‑masters?
\end{itemize}

\subsection*{Conclusion}

The hyper‑master is not ``an even more complicated version'' of the master series --- it is a \emph{necessary generalization} that establishes an exact link with the classical theory of hypergeometric functions. It is precisely in the form of the hyper‑master that the Master‑J formalism reveals its completeness, showing that ultra‑radicals, ultra‑logarithms and hypergeometric functions are particular manifestations of a single algebraic structure.

\section{Universal master series inside and outside the radius of convergence}
\label{sec:third_millennium}

Construction of power series for specific examples revealed that beyond the radius of convergence there exist not one, but \textit{two} distinct ways of parameter substitution into the same series.

\subsection{Canonical form and master series}
For a trinomial $Ay^d + By^g + C = 0$, we introduce the canonical form
\begin{equation}
p y^a = q + x y^b, \label{eq:canonical_trinomial}
\end{equation}
where $a, b \neq 0$, $a \neq b$. The solution is sought in the form
\begin{equation}
y = v \cdot M(1; a; b; z), \label{eq:ultraradical_sol}
\end{equation}
where the quantity $v$ is defined by
\begin{equation}
v^a = \frac{q}{p}, \label{eq:v_def}
\end{equation}
and the series argument is given by
\begin{equation}
z = \frac{x v^b}{a q}. \label{eq:z_def}
\end{equation}
The \textit{universal master series}, obtained via the Lagrange inversion formula for Eq.~\eqref{eq:canonical_trinomial}, can be written as a sum whose coefficients are expressed through a product of factors:
\begin{equation}
M(m; a; b; z) = m + z + \frac{(m-a+2b)}{2!}z^2 + \frac{(m-a+3b)(m-2a+3b)}{3!}z^3 + \dots \label{eq:master_series}
\end{equation}

The master series itself provides a unique solution within the circle of convergence. However, the quantities $v$ and $z$ defined by Eqs.~\eqref{eq:v_def} and~\eqref{eq:z_def} are multivalued functions of the exponents $a$ and $b$. Depending on which specific branches of the roots $v$ and $z$ are substituted into the master series, the corresponding root of the original trinomial equation is obtained. A key role in branch management is played by the auxiliary complex quantity
\begin{equation}
f = \frac{\ln|q/p| + i\bigl[\arg(q/p) + 2\pi N\bigr]}{a}, \label{eq:f_def}
\end{equation}
which depends on the chosen branch $N \in \mathbb{Z}$ of the complex logarithm. The base roots are expressed through it as $v = e^{f}$ and the series argument as $z = x e^{b f}/(a q)$. The quantity $f$ also uniquely determines the index of the sheet of the Riemann surface on which the sought root of the original equation is located. This index is computed by the formula
\begin{equation}
u = \left\lceil \frac{\operatorname{Im}(f)}{2\pi} - \frac{1}{2} \right\rceil, \label{eq:u_def}
\end{equation}
where $\lceil \cdot \rceil$ denotes rounding up to the nearest integer.

Equations with fractional, irrational, or complex exponents possess an important feature: many of their formal roots do not satisfy the original equation if exponentiation to powers $a$ and $b$ is performed via the principal branch of the complex logarithm during verification. Consider, for example, the equation
\begin{equation}
y^{2/3} + 0.01\, y^{1/2} + 1 = 0. \label{eq:example_trinomial}
\end{equation}
Under standard verification via the principal branch, none of the obtained roots reduces the equation to an identity. However, each root is strictly valid if the correct branch $u$ is used during substitution. The verification condition (or construction of an initial approximation for iterative refinement) is written as
\begin{equation}
\exp\!\left(\frac{2}{3}\bigl[\ln|y| + i(\arg y + 2\pi u)\bigr]\right) + 0.01\,\exp\!\left(\frac{1}{2}\bigl[\ln|y| + i(\arg y + 2\pi u)\bigr]\right) + 1 = 0, \label{eq:verification}
\end{equation}
where the parameter $u$ is applied identically to both power-law terms, despite the fact that only the exponent $a$ appears in definition~\eqref{eq:f_def}.

The roots $v$~\eqref{eq:v_def} act as base points on the Riemann surface: they fix the sheet and localize regions where solutions are located. At the same time, the fractional nature of the exponent $b$ determines the branching multiplicity: the number of distinct values of the argument $z$ generated by each base point $v$ is given by the denominator of the fraction $b$, and each such value generates an independent branch of the solution to the trinomial~\eqref{eq:canonical_trinomial}.

\subsection{Three canonical transformations}
The original trinomial is reduced to canonical form in three ways, corresponding to different dominant balances of terms:
\begin{align}
\text{AB:} \quad & Ay^d = -C - By^g \quad \Rightarrow \quad p=A,\; q=-C,\; x=-B,\; a=d,\; b=g, \nonumber \\
\text{BC:} \quad & By^{g-d} = -A - Cy^{-d} \quad \Rightarrow \quad p=B,\; q=-A,\; x=-C,\; a=g-d,\; b=-d, \nonumber \\
\text{CA:} \quad & Cy^{-g} = -B - Ay^{d-g} \quad \Rightarrow \quad p=C,\; q=-B,\; x=-A,\; a=-g,\; b=d-g. \label{eq:three_transforms}
\end{align}
Each method yields a different number of branches, since the degree of the equation for $v$ is determined by the modulus of the exponent $a$ in the relation $v^a=q/p$. For example, for $Ay^7 + By^4 + C = 0$, method AB gives all 7 roots, BC gives 3 roots, and CA gives 4 roots. With the argument of $z$ held fixed but its modulus smoothly increasing, the independent argument may exit the radius of convergence of method AB, which requires switching to the alternative transformation BC or CA; together they also yield all 7 roots.

\subsection{Geometric criterion for analytic continuation}

The key task is to determine which root from transformations BC or CA is the analytic continuation of the current branch $n$ of transformation AB. Despite the power of modern computational techniques, finding patterns amenable to simple formulas proved impossible without the aid of modern AI tools. A strict geometric criterion was formulated only in December 2025.

Mathematically, it is implemented in several steps.

\subsubsection{Reduction to unit canonical form}
For unambiguous numbering of branches $n$ in transformation AB, the original equation is first reduced to unit form. To this end:
\begin{enumerate}
    \item if necessary, the equation is multiplied by the unknown $y$ raised to such a power that the real parts of all exponents become positive;
    \item the coefficient of the term with the largest real part of the exponent is normalized to $1$;
    \item the constant term coefficient is normalized to $-1$.
\end{enumerate}
As a result, an equation of the form
\begin{equation}
y^d + B y^g - 1 = 0. \label{eq:canonical_unit}
\end{equation}
is obtained. Then transformation AB (see~\eqref{eq:three_transforms}) yields
\[
y^d = 1 - B y^g,
\]
whence $v^d = 1$, i.e., the base points are the $d$-th roots of unity.

\subsubsection{Sector partition of the Riemann surface}
If the exponents are real, the sector centers are given by the simplest formula:
\[
L(n) = \frac{2\pi b n}{a}.
\]
Each sector is bounded by the interval
\[
L(n) - \frac{\pi b}{a} < \arg z < L(n) + \frac{\pi b}{a}.
\]

For $|z| < R$ (inside the circle of convergence), sector partition is not required: computation proceeds directly via $f(n)$.

For $|z| \ge R$, it is necessary to determine which branch $h$ (transformation BC) or $k$ (transformation CA) is the analytic continuation of the given branch $n$.

\subsubsection{Alternative transformations}
From substitutions~\eqref{eq:three_transforms} we obtain:

\begin{align}
\text{BC:} \quad & B y^{g-d} = -1 + y^{-d}, \quad p = B,\; q = -1,\; x = 1,\; a = g-d,\; b = -d, \quad N = h, \nonumber \\
\text{CA:} \quad & -y^{-g} = -B - y^{d-g}, \quad p = -1,\; q = -B,\; x = -1,\; a = -g,\; b = d-g, \quad N = k. \nonumber
\end{align}

For transformation BC:
\[
v^{g-d} = -\frac{1}{B}, \qquad
f(h) = \frac{\ln|1/B| + i(\arg(-1/B) + 2\pi h)}{g-d}, \qquad
L(h) = \Im(bf(h)).
\]

For transformation CA:
\[
v^{-g} = B, \qquad
f(k) = \frac{\ln|B| + i(\arg(B) + 2\pi k)}{-g}, \qquad
L(k) = \Im(bf(k)).
\]

\subsubsection{Branch selection criterion}

The frequency of points $n$ equals the sum of frequencies of points $h$ and $k$:
\[
\frac{a}{2\pi} = \frac{a-b}{2\pi} + \frac{b}{2\pi}.
\]
Consequently, within a single sector $n$, points $L(h)$ and $L(k)$ cannot simultaneously reside. Therefore, analytic continuation is determined uniquely: the branch $h$ or $k$ that falls inside sector $n$ is the continuation of branch $n$ of transformation AB.

\subsubsection{Special case: sector boundary}

On the boundary between two sectors $n$ and $n+1$, points $L(h)$ and $L(k)$ may simultaneously appear (e.g., $h=0$ and $k=0$ on the boundary between $n=0$ and $n=-1$). In this case, the geometric criterion does not ensure a unique choice, since both branches intersect on the radius of convergence.

To synchronize numbering with the Lambert $W$-function, an additional rule is introduced:
\begin{quote}
\textit{If on the boundary of sectors $n$ and $n+1$ branches $h$ and $k$ are simultaneously present, then the branch with index $n+1$ is continued via the $h$-series, and the branch with index $n$ via the $k$-series.} A rigorous proof of this convention and its general statement are given in~\cite{Berezin2026_ULog} (Theorem~1).
\end{quote}

This rule ensures synchronization of ultraradical branch numbering with branches $W_0$ and $W_{-1}$ of the Lambert function, which is critically important for problems where a physical process may be described either via the ultraradical or via the $W$-function (e.g., adiabatic and isothermal limits in thermodynamics).

\subsection{How to use the ready-made function}

The entire algorithm described above --- selection of canonical transformation, analytic continuation, switching between series, and geometric branch criterion --- is implemented in program code and does not require deep expertise from the user. To apply the ultraradical to one's problem, it suffices to perform a few simple steps.

\subsubsection{Reduction to unit canonical form}

The original equation must be reduced to the canonical form
\begin{equation}
y^a = 1 + a x y^b, \label{eq:user_canonical}
\end{equation}
with the mandatory condition
\begin{equation}
a > b > 0. \label{eq:condition_ab}
\end{equation}

If the equation coefficients have a different form, they can be transformed via the substitution:
\[
w = y \left( \frac{q}{p} \right)^{\frac{1}{a}},\qquad
z = a q x \left( \frac{p}{q} \right)^{\frac{b}{a}},
\]
which maps the equation
\[
p w^{a} = q + z w^{b}
\]
into canonical form~\eqref{eq:user_canonical}.

The condition \( a > b > 0 \) is not an artificial restriction. It is necessary for correct operation of the geometric branch selection criterion: it is precisely the parameter \( a \) that determines the partition of the imaginary axis into sectors, and violation of this inequality distorts the branching structure.

\subsubsection{Function call}

After reduction to canonical form, the ultraradical function is called with three parameters:
\[
y = \sqrt[{n; a; b}]{x},
\]
where:
\begin{itemize}
    \item \( a \) and \( b \) are the exponents from the canonical equation;
    \item \( x \) is the argument (from the same equation);
    \item \( n \in \mathbb{Z} \) is the branch index (\( n = 0 \) for the principal branch).
\end{itemize}

All transformations related to the choice of analytic continuation method, switching between series, and branch synchronization occur automatically inside the function. The user does not need to know which series (AB, BC, or CA) is used at a given moment --- this is decided by the program code based on the geometric criterion.

\section{Additional Parameter m Possibilities}
\label{app:parameter_m}

Parameter $m$ can take values beyond 0 and 1. Using the power identity:
\[
M^{m}(1;a;b;x) = M(1;a/m;b/m;mx)
\]
we obtain power series allowing different interpretation of degree $m/r$.

\begin{multline*}
M^{m}(1;a;b;x) = M(1;a/m;b/m;mx) = \\
1 + m\left(x + \frac{x^{2}}{2}(m - a + 2b) + \frac{x^{3}}{3!}(m - a + 3b)(m - 2a + 3b) + \cdots\right)
\end{multline*}

Example:
\[
(1 + ax)^{\frac{m}{a}} = M^{m}(1;a;0;x) = M(1;a/m;0;mx) = 1 + m\left(x + \frac{x^{2}}{2}(m - a) + \frac{x^{3}}{3!}(m - a)(m - 2a) + \cdots\right)
\]

Numerical example:
\[
(1 + 3 \cdot 0.01)^{\frac{2}{3}} = M^{2}(1;3;0;0.01) = 1 + 2\left(0.01 + \frac{0.01^{2}}{2}(2 - 3) + \frac{0.01^{3}}{3!}(2 - 3)(2 - 6) + \cdots\right)
\]
\[
= 1 + 0.02 - 0.0001 + \frac{0.000004}{3} - \frac{0.00000007}{3} + \cdots
\]

\subsection{The Master Core}

To resolve ambiguity associated with the constant term, we introduce the concept of the \textbf{master core} -- the power series without the constant term:

\[
\text{Core}(m; a; b; x) = \sum_{\ell=1}^{\infty} \frac{x^\ell}{\ell !} \prod_{\gamma=1}^{\ell -1} (m - a\gamma + b\ell)
\]

The complete master series is then expressed as:

\[
M(m; a; b; x) = m + \text{Core}(m; a; b; x)
\]

For $m = 0$, the core coincides with the full series; for $m = 1$, it describes "oscillations" around unity. This separation simplifies the analysis and transformation of master series, providing greater flexibility in both theoretical and applied contexts.

\subsection{Super-Master Series: Generalized Parameter c}

To describe a broader class of problems, the concept of a \textbf{super-master series} is introduced, which depends on an additional generalized parameter $c$. This definition naturally extends the basic master series:

\begin{align}
S(m;a;b;x;c) &= m + c\bigg( x + \sum_{\ell=2}^{\infty} \frac{x^{\ell}}{\ell!} \prod_{\gamma=1}^{\ell-1} (cm - a\gamma + b\ell) \bigg) \notag \\
&= m + c\bigg( x + \frac{x^{2}}{2}(cm - a + 2b) \notag \\
&\quad + \frac{x^{3}}{3!}(cm - a + 3b)(cm - 2a + 3b) \notag \\
&\quad + \frac{x^{4}}{4!}(cm - a + 4b)(cm - 2a + 4b)(cm - 3a + 4b) + \cdots \bigg).
\label{eq:super_master_general}
\end{align}

In particular, for canonical cases, we obtain:
\begin{align}
M^{c}(1;a;b;x) &= S(1;a;b;x;c) = 1 + c\bigg( x + \frac{x^{2}}{2}(c - a + 2b) \notag \\
&\quad + \frac{x^{3}}{3!}(c - a + 3b)(c - 2a + 3b) + \cdots \bigg),
\label{eq:super_master_c1}
\end{align}
\begin{align}
cM(0;a;b;x) &= S(0;a;b;x;c) = c\bigg( x + \frac{x^{2}}{2}(2b - a) \notag \\
&\quad + \frac{x^{3}}{3!}(3b - a)(3b - 2a) + \cdots \bigg).
\label{eq:super_master_c0}
\end{align}

It should be noted that this formalism possesses a high degree of generality and unification. In particular, the results of D\v{z}. Belki\'c \cite{Belkic2018}, who obtained compact analytical expressions for arbitrary powers of trinomial equation roots via the confluent Fox--Wright function (see formula 11.15 in the cited work), are naturally derived from the super-master series as a special case for a fixed value of the leading degree parameter. This highlights the fundamental nature of the proposed analytical apparatus.
\paragraph{Connection with the series of Lambert and Euler.}
As shown in the review by Corless et al.~\cite{Corless1996}, Euler in 1779 obtained a series for the solution of the equation $y^A - y^B = (A-B)xy^{A+B}$ in the form
\[
y^c = 1 + cv + \frac{c(c+A+B)}{2!}v^2 + \frac{(c+A+2B)(c+2A+B)}{3!}v^3 + \cdots
\]
(formula 1.2 in~\cite{Corless1996}). Dividing Euler's equation by $y^B$ and making the substitution $a = A-B$, $b = A$ (so that $A = b$, $B = b-a$), we obtain the canonical master equation $y^a = 1 + axy^b$, and Euler's series for $y^c$ transforms into the super-master series $S(1; a; b; x; c)$.

Thus, in historical retrospect:
\begin{itemize}
\item the \textbf{master series} $M(1; a; b; x)$ corresponds to the \emph{Lambert series} (1758) for the root of a trinomial equation;
\item the \textbf{super-master series} $S(1; a; b; x; c)$ corresponds to the \emph{Euler series} (1779) for an arbitrary power of the root.
\end{itemize}

\section{Software Implementation and Analytic Continuation Algorithm}
\label{app:implementation_and_algorithm}

\subsection{Overview of the Computational Framework}

The Master-J method has been implemented as a comprehensive computational framework in the Maple computer algebra system. The implementation provides both symbolic and numerical tools for working with master series and ultra-radicals, with particular emphasis on deterministic analytic continuation beyond the radius of convergence.

The core implementation is publicly available at two Zenodo repositories:
\begin{itemize}
    \item \textbf{SuperMaster (General Framework):} \href{https://doi.org/10.5281/zenodo.17717360}{\texttt{10.5281/zenodo.17717360}}
    \item \textbf{Ultra-Radical Algorithm:} \href{https://doi.org/10.5281/zenodo.17743595}{\texttt{10.5281/zenodo.17743595}}
\end{itemize}

\subsection{Algorithm for computing the master series}
The super-master series is computed using the recurrence formula
\begin{equation}
S = m + c x + \sum_{\ell=2}^{\ell_{\max}} \Delta_\ell, \quad
\Delta_\ell = \frac{c \cdot x^\ell}{\ell!} \prod_{\gamma=1}^{\ell-1} \left( c m - a \gamma + b \ell \right),
\label{eq:master_series_rec}
\end{equation}
summation terminates when the relative convergence criterion is satisfied
\begin{equation}
\left| \frac{\Delta_\ell}{S_{\ell-1}} \right| < \varepsilon,
\label{eq:conv_criterion}
\end{equation}
where $\varepsilon = 0.5 \cdot 10^{-D}$ sets the required precision of $D$ significant digits.
Upon reaching $\ell_{\max}$ or a zero increment (binomial expansion), the algorithm returns the current sum.

The presented recurrence scheme is universal for the generalized master series, where for $m=0$ the parameter $c$ sets a scale factor, and for $m=1$ the parameter $c$ sets the power. For the ultraradical case ($m=1$, $a \neq 0$, $b \neq 0$, $a \neq b$) the super-master number $N_\ell=\prod_{\gamma=1}^{\ell-1} \left( c m - a \gamma + b \ell \right)$ can be expressed in closed form via the gamma function:
\begin{equation}
N_\ell = (-a)^{\ell-1} \, 
\frac{\Gamma\!\left( \dfrac{(a-b)\ell - c m}{a} \right)}
     {\Gamma\!\left( \dfrac{a - c m - \ell b}{a} \right)},
\label{eq:M_L_gamma2}
\end{equation}
or in a numerically stable form:
\begin{equation}
N_\ell = a^{\ell-1} \, 
\frac{\Gamma\!\left( \dfrac{c m + \ell b}{a} \right)}
     {\Gamma\!\left( \dfrac{c m + \ell b}{a} - \ell + 1 \right)},
\label{eq:M_L_gamma_clean2}
\end{equation}
This representation allows analytical detection of zero series terms: if the argument of the gamma function in the denominator takes a non-positive integer value, the corresponding series term vanishes.

\paragraph{Auxiliary functions.}
\begin{itemize}
    \item \texttt{O9(x)} --- rounding with correct nine-carry propagation; ensures stable branch numbering under machine precision.
    \item \texttt{UltraRadical(n, a, b, x, M)} --- computation of the $n$-th branch of the ultraradical; automatically selects one of the three canonical parameter-substitution methods depending on the argument position relative to the radius of convergence.
    \item \texttt{URad(n, a, b, x, M)} --- hybrid root refinement: switches to Newton or Schr\"oder method when the master series converges slowly (near the radius of convergence).
\end{itemize}

The full implementation of the algorithm in Python, including the functions \texttt{O9}, \texttt{UltraRadical} and \texttt{URad}, as well as example calculations for the gravitational thermocompressor, is available in an open repository with a persistent DOI identifier:
\begin{equation}
\texttt{https://doi.org/10.5281/zenodo.20677790} .
\end{equation}

\section{Structural Patterns in Algebraic Equations}

Analysis of ultra-radical power series reveals the external structure of root formulas for equations $y^\alpha = q + xy^\beta$.

\subsection{Universal Decomposition Method}

For irreducible equations where $\gcd(\alpha,\beta)=1$, the solution admits a structural decomposition:
\[
y = f_0 + f_1 + f_2 + \cdots + f_{\alpha-1}
\]
where each term follows the pattern:
\[
f_j = L_j \cdot V_j^{(1-(\alpha-\beta)j)/\alpha}, \quad L_j = \left(\frac{x}{\alpha}\right)^j \cdot \frac{\prod_{K=1}^{j-1}(1 + \beta j - \alpha K)}{j!}
\]

\subsection{Detailed Example: Cubic Case y^3 = q + xy}

\subsubsection{Step 1: Structural Decomposition}
For $y^3 = q + xy$, we have $\alpha=3$, $\beta=1$, $t=\alpha-\beta=2$:
\[
y = f_0 + f_1 + f_2 = V_0 + \frac{x}{3V_1} + 0
\]
since $L_2=0$ for this case.

\subsubsection{Step 2: Root Representation}
Let $V_0 = a$, $V_1 = b$. For each root we apply cyclic symmetries:
\begin{align*}
y_0 &= a_0 + \frac{x}{3b_0} \\
y_1 &= a_1 + \frac{x}{3b_1} \\
y_2 &= a_2 + \frac{x}{3b_2}
\end{align*}
with $a_1 = \omega a_0$, $a_2 = \omega^2 a_0$, $b_1 = \omega b_0$, $b_2 = \omega^2 b_0$, where $\omega = e^{2\pi i/3}$.

\subsubsection{Step 3: Apply Vieta's Formulas}
\begin{align*}
\text{Sum: } & y_0 + y_1 + y_2 = (a_0 + a_1 + a_2) + \frac{x}{3}\left(\frac{1}{b_0} + \frac{1}{b_1} + \frac{1}{b_2}\right) = 0 \\
\text{Pairwise: } & y_0y_1 + y_0y_2 + y_1y_2 = -x \\
\text{Product: } & y_0y_1y_2 = q
\end{align*}

\subsubsection{Step 4: Solve the System}
From the sum condition and symmetry, we find $a_0 = b_0$ (thus $V_0 = V_1$). 

Substituting into the product condition:
\[
\left(a_0 + \frac{x}{3a_0}\right)\left(\omega a_0 + \frac{x}{3\omega a_0}\right)\left(\omega^2 a_0 + \frac{x}{3\omega^2 a_0}\right) = q
\]
Simplifying yields:
\[
a_0^6 - a_0^3q + \frac{x^3}{27} = 0
\]
Thus $V_0 = a_0$ is determined by this sextic equation.

\subsection{General Method}

For any equation $y^\alpha = q + xy^\beta$:
\begin{enumerate}
\item Write the structural decomposition with parameters $V_0, V_1, \ldots, V_{\beta-1}$
\item Express all roots using symmetry transformations
\item Apply Vieta's formulas to obtain equations
\item Solve the resulting system for $V_j$
\end{enumerate}

This constructive approach generates root formulas whose power series match the corresponding ultra-radicals, revealing the fundamental patterns underlying algebraic solutions.

\subsection{Quintic Case and Beyond}

For equations of degree 5 and higher, the same method applies but leads to fundamental limitations.

\subsubsection{Quintic Case: y^5 = q + 5xy}

The structural decomposition gives:
\[
y = V_0 + \frac{x}{V_1^3} - \frac{x^2}{V_2^7} + \frac{x^3}{V_3^{11}}
\]
with parameters $a=V_0$, $b=V_1$, $c=V_2$, $d=V_3=abc$.

Applying the method yields a system of equations. After substitutions $X=bc$, $Y=a^2b$, $Z=a^3c$, $P=-x^5$, we obtain:
\begin{align*}
&Y^2(XZ+P) - XZ^2 + PZ = 0 \\
&Y^4(PX) + Y^2(-X^3Z + X^2Z^2 + 3PXZ - P^2) - PXZ^2 = 0 \\
&Y^6(-P^2X) + Y^4(X^3Z^3 - 10PX^2Z^2 + 10P^2XZ) - Y^3(X^4Z^2q) \\
&\quad + Y^2(-10PX^2Z^3 - 10P^2XZ^2 - P^3Z) - PX^2Z^4 = 0
\end{align*}

\subsubsection{Fundamental Limitation}

This demonstrates why quintic equations cannot be generally solved in radicals: the structural approach necessarily leads to resolvents of higher degree than the original equation. The parameter system for degree 5 generates equations of degree up to 12, confirming the Abel-Ruffini theorem through constructive means.

The ultra-radical $\sqrt[n;5;1]{x}$ thus represents the minimal analytic continuation of this structural pattern beyond radical solvability.
\section{Open Problems and Future Directions}

The Master-J framework, as presented in this work, naturally gives rise to a set of fundamental questions that delineate its potential scope and place within modern mathematics. The following open problems are proposed to guide future research and discussion within the scientific community.

\subsection{Theoretical Foundations}

\begin{enumerate}
    \item \textbf{Universality of the Master Series}: To what extent can all elementary and special functions be represented as specific instances of the master series $M(m;s;r;x)$? Which functions, if any, resist such representation and why?
    
    \item \textbf{The Merge Operation}: Does the $@$ (merge) operation constitute a fundamental mathematical operation, analogous to addition or multiplication? Can a consistent algebra be constructed based on this operation?
    
    \item \textbf{Analytical Continuation}: To what extent can the geometric criterion for branch selection be extended to provide deterministic analytic continuation for other classes of multi-valued functions beyond their convergence radii while maintaining branch continuity?
\end{enumerate}

\subsection{Computational and Practical Scope}

\begin{enumerate}
    \setcounter{enumi}{3}
    \item \textbf{Limitations of the Method}: What are the fundamental limitations of the Master-J method? For which classes of equations or problems is it provably inapplicable?
    
    \item \textbf{Comparative Advantage}: For which computationally challenging problems (e.g., in quantum mechanics with complex exponents, chaotic systems, or nonlinear optics) does Master-J offer a tangible advantage over iterative methods, particularly regarding branch identification and the absence of initial guesses?
\end{enumerate}

Addressing these questions will not only validate the method presented but also chart the course for its evolution into a potential new paradigm within computational mathematics and analysis.

\section*{For the Latest Developments}

Readers interested in following the progress on these open problems and accessing the most recent version of this research are encouraged to visit the Zenodo repository using the Concept DOI below. This link always resolves to the latest version and may contain updates, corrections, and additional materials:

\begin{center}
\href{https://doi.org/10.5281/zenodo.17682133}{\large\texttt{10.5281/zenodo.17682133}} , ru: \href{https://doi.org/10.5281/zenodo.17822264}{\large\texttt{10.5281/zenodo.17822264}}
\href{https://doi.org/10.5281/zenodo.22995073}{\large\texttt{10.5281/zenodo.22995073}} , ru: \href{https://doi.org/10.5281/zenodo.22994577}{\large\texttt{10.5281/zenodo.22994577}}
\begin{itemize}
\item \textbf{Paper 1 (this work):} ultra-radical --- arXiv:2512.18643; Zenodo: \href{https://doi.org/10.5281/zenodo.22995073}{en}, \href{https://doi.org/10.5281/zenodo.22994577}{ru}.
\item \textbf{Paper 2:} master series and merging --- arXiv:2609.33459; Zenodo: \href{https://doi.org/10.5281/zenodo.22995072}{en}, \href{https://doi.org/10.5281/zenodo.22994576}{ru}.
\item \textbf{Paper 3:} ultralogarithm, rigorous criterion --- arXiv:2611.23317; Zenodo: \href{https://doi.org/10.5281/zenodo.24219280}{en}, \href{https://doi.org/10.5281/zenodo.24218691}{ru}.
\end{itemize}
\end{center}

\subsection{On the Systematic Nature of the Master‑J Framework}

The Master‑J method provides a parametrized series representation for solutions to certain algebraic and transcendental equations. This structured approach enables:
\begin{itemize}
    \item uniform convergence analysis across different equation types,
    \item a deterministic algorithm for analytic continuation,
    \item natural generalization to equations with arbitrary coefficients and multiple terms.
\end{itemize}
Such a unified viewpoint can streamline both theoretical analysis and practical implementation, much as the theory of hypergeometric functions unified solutions to linear differential equations.

\subsection{On the Role of Iterative Methods}

One might ask whether the Master-J framework renders iterative methods obsolete.
The answer is no—rather, it complements them. Each approach has its domain of excellence:

• \textbf{Master-J excels} when $|x| \ll R$ (far from the convergence boundary), for complex or fractional exponents, and when branch continuity or analytical insight is required.

• \textbf{Iterative methods remain valuable} for simple polynomials with small integer degrees, near the convergence radius, or in hybrid schemes where Master-J provides optimal initial guesses.

The emergence of new concepts and complex objects in mathematics does not replace the existing language—it \textbf{expands} it. The Master-J approach offers an alternative technique that can complement existing computational methods

\section*{Historical Note: The Search for the Ultra-Radical and the Role of Computational Tools}

The quest for an analytic solution to generalized trinomial equations has a long and distinguished history. The problem essentially revolves around inverting the function defined by $y^a = 1 + a x y^b$, a quest for what we now term the \emph{ultra-radical}.

The journey began in the 18th century with the pioneering work of **Johann Heinrich Lambert** (1758) \cite{Mezo2022}, who systematically studied the trinomial $x^m + p x = q$. He derived successive approximations and recognized the structure of a series solution, laying the groundwork for later investigations. **Leonhard Euler**, inspired by Lambert, examined the more general equation $x^\alpha - x^\beta = (\alpha-\beta) v x^{\alpha+\beta}$ and its limiting logarithmic form. These efforts highlighted the inherent complexity of the problem but stopped at deriving series expansions for specific instances, lacking a unified, parametrized formalism and a method for analytic continuation.

In the late 18th century, **Erland Samuel Bring**, in his work on the quintic equation, discovered a crucial reduction: the general quintic can be transformed into the form $y^5 + p y + q = 0$. This is equivalent to the ultra-radical equation with parameters $a=5, b=1$. Thus, Bring's radical, a key element in the theory of quintic solvability, is a distinct, isolated point $(5,1)$ in the two-dimensional parametric space of the ultra-radical.

For centuries, these contributions remained as separate, insightful fragments. A complete, unified theory required not only conceptual leaps but also a tool capable of handling the immense symbolic complexity involved. The manual derivation of the general series coefficients, the identification of the three isoroot transformations (AB, BC, CA), and the development of the merge operation (@) for multi-term equations represented a combinatorial and algebraic challenge of staggering proportions.

The modern resolution of this centuries-old search became possible only with the advent of powerful computational tools. The author's investigation illustrates this point dramatically:
\begin{itemize}
    \item The complete parametric power series and the three fundamental transformations were derived and verified using symbolic computation within one week in 2018.
    \item In contrast, the **merge operation (@)**, a generalization essential for solving equations with an arbitrary number of terms, eluded manual derivation for five years. It was finally identified and formulated in 2023, again through an iterative cycle of computational experimentation and pattern recognition aided by the structural insights provided by the system.
    \item In contrast, the \textbf{geometric criterion for conjugate series selection} --- the key to deterministic analytic continuation beyond the convergence radius --- eluded discovery for \textbf{seven years}. Despite extensive numerical experiments with various branch parametrizations, no closed-form pattern emerged. The breakthrough came only in December 2025, when the sector structure of the Riemann surface was finally recognized as the organizing principle. This criterion is what transforms the master series from a local tool into a global analytical function.
    \item The recent preparation of this manuscript, including the verification of complex identities, branch continuity proofs, and structuring of the theory, was significantly accelerated by using large language models (LLMs) specialized in mathematical reasoning. These tools acted as tireless assistants in code generation, consistency checking, and suggesting expositional structures.
\end{itemize}

Historical analysis reveals that the ultraradical is a \emph{computationally native} object. Its rigorous description remained inaccessible to Lambert, Euler, and Bring not due to a lack of mathematical intuition, but owing to the computational limitations of their era. The present work realizes their vision by leveraging modern symbolic and numerical tools.

\section*{Acknowledgements}

The author expresses sincere gratitude to:

- \textbf{Professor Alexander E. Dubinov} (Russian Federal Nuclear Center–VNIIEF and Sarov Institute of Physics and Technology) for stimulating discussions and for drawing attention to the potential applications of the ultra-radical in plasma physics.

- \textbf{Pavel V. Berezin, Alexey V. Berezin, and Andrey V. Gruzdov} for their extensive work on computations, numerical verification, and fruitful discussions.

\textbf{Professor Tuyen Truong} (University of Oslo, Norway) and \textbf{Professor Abdulaziz Alhaidari} (Saudi Center for Theoretical Physics) for their valuable feedback, critical questions, and constructive discussions that greatly improved this work.

\paragraph*{Recommendation on AI Tool Usage.}
During the preparation and editing of this manuscript, the author actively utilized the AI language models DeepSeek Math and Qwen. These tools played a significant role in structuring complex mathematical derivations, verifying formula consistency, optimizing algorithmic descriptions, and refining the academic quality of the text. Researchers in mathematical analysis, special functions, and computational mathematics are encouraged to leverage such systems to accelerate analytical investigations, cross-check proofs, and enhance the clarity and rigor of scientific publications.
Analytical computations and numerical verification were performed using the computer algebra systems \textbf{Wolfram Alpha} and \textbf{Maple}, as well as the \textbf{Python} programming language with its scientific libraries.

\vspace{1cm}
\noindent
\textbf{Document version:} 4.0
\noindent\begin{itemize}
\item \textbf{Notation system change:} In versions 1.x, the parameter set $(m; s; r; x)$ was used, where $r=a$ and $s=b$ played roles analogous to linear coefficients. Starting from version 2.0, for convenience when working with algebraic equation exponents, the system $(m; a; b; x)$ has been adopted, where:
  \begin{itemize}
  \item $a$ — the first exponent (main degree in $y^a$)
  \item $b$ — the second exponent (in the $y^b$ term)
  \item This alignment makes the notation intuitively clear when solving equations of the form $y^a = 1 + axy^b$ or $AY^a + BY^b + C = 0$
  \end{itemize}

\textbf{Current implementation:} \href{https://doi.org/10.5281/zenodo.17682133}{10.5281/zenodo.17682133} , ru: \href{https://doi.org/10.5281/zenodo.17822264}{10.5281/zenodo.17822264}
\end{itemize}

\end{document}